\documentclass[11pt]{article}
\usepackage{amssymb,amsmath,amsthm}
\usepackage{epsfig}
\usepackage{verbatim}
\usepackage{graphicx}

\usepackage{color}

\usepackage{MnSymbol}

\newtheorem{thm}{Theorem}[section]
\newtheorem{cor}[thm]{Corollary}
\newtheorem{prop}[thm]{Proposition}

\newtheorem{rem}[thm]{Remark}

\newtheorem{lemma}[thm]{Lemma}

\newcommand{\R}{\mathbb{R}}

\newcommand{\pv}{{\rm pv}}

\newcommand{\grad}{\nabla}

\newcommand{\dt}{\frac{d}{dt}}

\newcommand{\dxu}{\partial_{x_1}}

\newcommand{\te}{\theta}
\newcommand{\eps}{\epsilon}

\newcommand{\les}{\Lambda^{\gamma+s}}

\newcommand{\heatinv}{(\partial_t\!-\!\Delta)_0^{-1}}

\begin{document}

\title{Global regularity for 2D Boussinesq temperature patches \\ with no diffusion}

\author{Francisco Gancedo and Eduardo Garc\'ia-Ju\'arez}

\date{\today}

\maketitle

\begin{abstract}
This paper considers the temperature patch problem for the incompressible Boussinesq system with no diffusion and viscosity in the whole space $\R^2$. We prove that for initial patches with $W^{2,\infty}$ boundary the curvature remains bounded for all time. The proof explores new cancellations that allow us to bound $\grad^2u$, even for those components given by time dependent singular integrals with kernels with nonzero mean on circles.
In addition, we give a different proof of the $C^{1+\gamma}$ regularity result in \cite{DANCHINZHANG}, $0<\gamma<1$, using the scale of Sobolev spaces for the velocity. Furthermore, taking advantage of the new cancellations, we  go beyond to show the persistence of regularity for $C^{2+\gamma}$ patches.
\end{abstract}


\setcounter{tocdepth}{1}

\section{Introduction}

In this paper we consider the following active scalar equation 

\begin{equation}\label{temperature}
\te_t + u\cdot\grad\te =0,
\end{equation}
for incompressible fluids
\begin{equation}\label{incompressible}
\grad\cdot u=0,
\end{equation}
which depending on the physical context can be seen as the mass or energy conservation. In this latter case $\te=\te(x,t)$ corresponds to the temperature transported without diffusion by the fluid, which moves with velocity $u=(u_1(x,t),u_2(x,t))$.  We close the system with the Boussinesq model for the momentum equation
\begin{equation}\label{Boussinesq}
u_t +u\cdot\grad u-\nu\Delta u + \nabla P=g(0,\te),
\end{equation}
where $P$ is the pressure, $\nu>0$ the viscosity, $g$ the gravity, $x=(x_1,x_2)\in\R^2$ and $t\geq 0$. By a  change of variables we simplify matters by taking gravity $g = 1$ and viscosity $\nu = 1$. 

\par
This system first arose as a model to study natural convection phenomena in geophysics \cite{PEDLOSKY}, \cite{MAJDA}, as for example in the very important Rayleigh-B\'enard problem \cite{CONSTANTINDOERING}. There, the density variations can frequently be neglected except in the buoyancy term, avoiding in this way the calculation of sound waves.  From the mathematical point of view, the interest resides on its connection to the Navier-Stokes and Euler equations since it presents  vortex stretching \cite{MABE}.
\par
For that reason, the well-posedness of this model has recently attracted a lot of attention, starting with the results of Chae \cite{CHAE} and Hou and Li \cite{HOULI} for regular initial data in the whole space $\R^2$. Later, making use of paradifferential calculus techniques, results for rougher initial data appeared. In particular, Abidi and Hmidi \cite{ABIDIHMIDI} established the global well-posedness for initial data in the Besov space $B_{2,1}^0$ (see Appendix for definition), then Hmidi and Keraani \cite{HMIDIKERAANI} proved global existence and regularity for initial data $\theta_0\in L^2$, $u_0\in H^s, s\in[0,2)$ and finally Danchin and Paicu obtained the uniqueness \cite{DANCHINPAICU}. The persistence of regularity in Sobolev spaces was completed by Hu, Kukavica and Ziane in \cite{KUKA}. Available global-in-time results in three dimensions require the initial data to be small \cite{DANCHINPAICU}, as the system includes Navier-Stokes equations as a particular case. 
\par
The Boussinesq system is also used to model large scale atmospheric and oceanic flows, where the viscosity and diffusivity constants are usually different in the horizontal and vertical directions.
In these situations there are similar results for the case with anisotropic dissipation \cite{WU1}, \cite{WU2}, \cite{TITI}, \cite{LITITI}, and positive diffusivity but no viscosity \cite{CHAE}, \cite{HMIDIKERAANI2}, including results with Yudovich type initial data \cite{DANCHINPAICU2}. 
In contrast, the global well-posedness of the full inviscid case remains still as an open question, mathematically analogous indeed to the incompressible axi-symmetric swirling three-dimensional Euler equations \cite{MABE}. Simulations first indicated the possibility of finite-time blow-up, but there were also numerical evidence in the opposite direction \cite{WEINAN}. Recently, based on numerical studies, a new scenario for finite time blow-up in 3D Euler equations has been proposed \cite{LUOHOU}. This situation has been adapted to rigorously prove the existence of finite time blow-up first for a 1D model of the Boussinesq equations \cite{KISELEV2} and  more recently for a modified version of the two dimensional case including incompressibility  \cite{KISELEV}. 
\par
Related to these problems, we consider here a case with singular initial data: the so-called {\em{Boussinesq temperature patch problem}} for \eqref{temperature}-\eqref{Boussinesq}. 
From \eqref{temperature} and the definition of the particle trajectories,
\begin{equation*}
\left\{\begin{aligned}
\frac{dX}{dt}(x,t)&=u(X(x,t),t),\\
X(x,0)&=x,
\end{aligned}\right.
\end{equation*}
an initial temperature patch $\theta_0=1_{D_0}$, i.e., the characteristic function of a simply connected bounded domain $D_0$, will remain as a patch $\theta(t)=1_{D(t)}$, where $D(t)=X(D_0,t)$. Therefore, one may wonder whether its boundary preserves the initial regularity.
This kind of problems were studied in the 80s for the well-known vortex patch problem. In that case, it was first thought from numerical results that there was finite-time blow-up, but later the global regularity was proved by Chemin \cite{CHEMIN} using paradifferential calculus and striated regularity techniques. The same result was proved by Bertozzi and Constantin \cite{CONST} in a geometrical way making use of harmonic analysis tools.
\par
Other counter dynamics scenarios to look for singularity formation come from the evolution of the interface between fluids of different characteristics. Used to model physically important problems such as water waves, porous media, inhomogeneous flow or frontogenesis, this contour dynamics setting has attracted a lot of attention in the last years. The appearance of finite time singularities was first proved for the Muskat problem \cite{CASTRO}, \cite{CASTRO2}, Euler \cite{CASTROEULER} and Navier-Stokes equations \cite{CASTRONAVIER} and from there different scenarios and results appeared \cite{COUTAND2}, \cite{COUTAND}.
For the SQG active scalar incompressible system \cite{CONSTMAJTAB} there is numerical evidence of pointwise collapse with curvature blow-up for the patch problem \cite{CORDOBANUM}. In addition,   it has been shown that the control of the curvature removes the possibility of pointwise interface collapse \cite{GANCEDO}.
\par
Recently Danchin and Zhang proved \cite{DANCHINZHANG} that if the initial temperature of the Boussinesq system \eqref{temperature}-\eqref{Boussinesq} is a $C^{1+\gamma}$ patch, it remains with the same regularity for all time. In this sense, we will denote $\partial D(t)\in C^{1+\gamma}$ if there exists a parametrization of the boundary 
\begin{equation}\label{parametrization}
\partial D(t) = \left\{z(\alpha,t)\in\R^2, \alpha\in[0,1]\right\}
\end{equation}
with $z(t) \in C^{1+\alpha}$.  
\par
 From the above, one may wonder if the curvature of a temperature patch in the Boussinesq system can blow up without self-intersection or if, on the contrary, it remains bounded for all time. We prove here that the latter occurs, that is, we show the persistence of $W^{2,\infty}$ regularity for Boussinesq temperature patches.

\medskip
\par
In \cite{DANCHINZHANG}, the authors use Besov spaces to measure regularity and the techniques of striated regularity to get the $C^{1+\gamma}$ propagation. The main idea is that to control the H\"older regularity of the patch one just needs to control the gradient of the velocity in certain directions, which can be translated into the vorticity equation and then treat this as a forced heat equation, hence achieving the gain of two derivatives. This result was done for a general $\te_0\in B_{q,1}^{2/q-1}$, $q\in(1,2)$ and then applied to temperature patches, as a patch always belongs to that space. 
\par
In this paper we exploit the fact of $\te$ being a patch. Indeed, one can get the H\"older persistence of regularity by controlling the $L^1(0,T;C^\gamma)$ norm of the gradient of the velocity using the particle trajectories of the system. This can be seen from the vorticity equation rewritten as
\begin{equation}\label{vortsplit}
\omega(t)=\omega_1+\omega_2+\omega_3:=e^{t\Delta}\omega_0-\heatinv \nabla\cdot(u\omega)+\heatinv\partial_1\te ,
\end{equation}
where $\heatinv f$ denotes the solution of the heat equation with force $f$ and zero initial condition:
\begin{equation*}\label{inverseheatequation}
(\partial_t-\Delta)^{-1}_0f:=\int_0^t e^{(t-\tau)\Delta}f(\tau)d\tau.
\end{equation*}
Above we use the standard notation $e^{t\Delta}f=\mathcal{F}^{-1}(e^{-t|\xi|^2}\hat{f})$, where $\hat{}$ and $\mathcal{F}^{-1}$ denote Fourier transform and its inverse.
\par
We note then that one could get $u\in L^1(0,T; B_{\infty,\infty}^2)$ by choosing more regular initial conditions, since the main limitation comes from the temperature term and singular integrals are bounded on Besov spaces. However, to control the boundedness of the curvature of the patch it will be needed to control the $L^\infty$ norm of the second derivatives of the velocity,
\begin{equation*}
\partial_k\partial_ju_i(t)=\partial_i^\perp \partial_j(-\Delta)^{-1}\left(e^{t\Delta}\partial_k\omega_0-\heatinv\partial_k \nabla\cdot(u\omega)+\heatinv\partial_k\partial_1\te  \right),
\end{equation*}
$(\partial_1^\perp,\partial_2^\perp)=(-\partial_2,\partial_1)$, i.e., we will need $u\in L^1(0,T;W^{2,\infty})$. This is not trivial since neither the Riesz transforms nor the operators $\partial_i\partial_j \heatinv$  are bounded for a general function in $L^\infty$.
In fact, it is known that this operator takes bounded functions to $BMO$  \cite{SCHLAG} (defined using parabolic cylinders instead of Euclidean balls). In addition, while one may expect to get rid of the Riesz transform by using striated regularity techniques (assuming that one could first get further regularity in the tangential direction to interpolate the $L^\infty$ norm), it would still be  necessary to bound $\Delta \heatinv \te$ (see Remark \ref{remstriated}). In fact, the associated kernel
\begin{equation}\label{laplacevort}
\frac{1}{4\pi t^2}\left(\frac{|x|^2}{4t}-1\right)e^{-|x|^2/4t}
\end{equation}
 is not integrable due to the singularity at the origin along parabolas $t=c |x|^2$. However, the cancellation due to the sign change through the parabola $t=|x|^2/4$ allows to appropriately define these operators as principal values. What is more, although this kernel has nonzero mean, we will show that it is bounded for $\te$ a patch. First we give a $C^{1+\gamma}$ result. Later we refine the idea used in \cite{CONST} to bound the gradient of the velocity, i.e., the combined fact that for a $C^{1+\gamma}$ patch the intersection of a small circle with its boundary is {\em{almost}} a semicircle and that the kernel was even with zero mean on circles. Although the latter is not true in our case (see kernel \eqref{laplacevort}), we encounter that the kernels present certain time-space cancellations on circles. In this scenario, the kernels now depend also on time so that the picture is no more static and therefore we need to take care of the evolution of the distance of the point to the boundary. 

The above result would prove that for $\te$ a patch $\grad \omega$ is bounded. We can polish the idea and adjust it to the operators $$\partial_i^\perp \partial_j(-\Delta)^{-1}\partial_k\partial_1\heatinv.$$

\medskip
The strategies above allow us to control the particle trajectories by estimating all the second derivatives of the velocity. There is hope that one can take advantage of the cancellation $W\cdot\grad\te\equiv 0$, where $W$ is a vector field tangent to the patch \cite{DANCHINZHANG}. Going beyond, we can take advantage of this cancellation and the $W^{2,\infty}$ result to prove the persistence of regularity for $C^{2+\gamma}$ patches. 

\medskip
In relation to the initial conditions, for $\omega_0\in H^s$ it is not possible to proceed as before for $s=0$ due to the term $\omega_1$ in \eqref{vortsplit}: specifically, one encounters the failure of the embedding $H^1(\R^2)\nhookrightarrow L^\infty$ and the restriction from the maximum gain of the heat equation shown in \cite{FERO2}. Nevertheless, if we don't restrict ourselves to Sobolev spaces, the result of persistence of regularity for the patch is still true for more general initial conditions (see Section \ref{Corollaries}). In this sense, the initial conditions to achieve the $C^{1+\gamma}$ result are at the same level of regularity to those in \cite{DANCHINZHANG} but in the Sobolev space scale.

\medskip
The structure of the paper is as follows: In the next section, we give a preliminary result of the $C^{1+\gamma}$ regularity using particle trajectories and Sobolev spaces for the velocity. Then, in Section \ref{mainthe} we present the main result: the control of the curvature. First we define the operators involved, then show the new cancellations encountered and finally prove that $u\in L^1(0,T;W^{2,\infty})$.  In Section \ref{Case33} we use the extra regularity of the patch in the tangential direction to show the persistence of $C^{2+\gamma}$ regularity. Additionally, in Section \ref{Corollaries} we give results for initial conditions in different spaces. In Appendix we include some results related to the heat and transport equations.

\section{Persistence of $C^{1+\gamma}$ regularity}\label{sec:2}

In this section the persistence of $C^{1+\gamma}$ regularity is proved. We state the result for velocity fields belonging to the Sobolev space $H^{\gamma+s}$, $s\in (0,1-\gamma)$ (see Section \ref{Corollaries} for $u_0\in B_{\infty,\infty}^{\gamma-1+s_1}\cap H^{s_2}$, $s_1\in(0,1-\gamma)$, $s_2\in(0,1)$). This is roughly at the same level of regularity needed for the result in \cite{DANCHINZHANG}, where for a general initial vorticity the striated regularity condition can be translated into $u_0\in B_{\infty, \infty}^{\gamma-1}\cap B_{q,1}^{-1+2/q}\hookrightarrow  B_{\infty, \infty}^{\gamma-1}\cap B_{2,1}^0$, with $q\in(1,\frac{2}{2-\gamma})$. 

\vspace{0.3cm}
\begin{thm}\label{Case1}
Assume $\gamma\in(0,1)$, $s\in(0,1-\gamma)$. Let $D_0\subset \R^2$ be a bounded simply connected domain with boundary $\partial D_0\in C^{1+\gamma}$, $u_0\in H^{\gamma+s}$ a divergence-free vector field and $\te_0=1_{D_0}$ the characteristic function of $D_0$. Then, there is a unique global solution $(u,\te)$ of \eqref{temperature}-\eqref{Boussinesq} such that
$$\te(x,t)=1_{D(t)}(x) \hspace{0.2cm}{\rm{and}} \hspace{0.2cm} \partial D\in L^\infty(0,T;C^{1+\gamma}),$$
where $D(t)=X(D_0,t)$ with $X$ the particle trajectories associated to the velocity field.\newline
Moreover, 
$$
u\in L^\infty(0,T;H^{\gamma+s})\cap L^2(0,T;H^{1+\gamma+s})\cap L^1(0,T; H^{2+\mu})\cap L^1(0,T;C^{1+\gamma+\tilde{s}}),
$$
for any $T>0$, $\mu<\min\{\frac12, \gamma+s\}$, $0<\tilde{s}<s$.  
\end{thm}

Proof: The first part of the proof consists of {\em{a priori}} estimates.  
From the transport equation for the temperature one gets
\begin{equation*}
\|\te(t)\|_{L^2}\leq \|\te_0\|_{L^2},\hspace{1cm}\|\te(t)\|_{L^\infty}\leq \|\te_0\|_{L^\infty}.
\end{equation*}
The basic energy inequality also holds
\begin{equation}\label{energy}
\frac12\dt\|u\|_{L^2}^2+\|\grad u\|_{L^2}^2\leq \|\te_0\|_{L^2}\|u\|_{L^2},
\end{equation}
so that by Gr\"onwall's lemma
\begin{equation}\label{uest}
\|u(t)\|_{L^2}^2\leq \left(\|u_0\|_{L^2}^2+\|\te_0\|_{L^2}^2\right)e^t-\|\te_0\|_{L^2}^2,
\end{equation}
\begin{equation}\label{graduest}
\int_0^t\|\grad u(\tau)\|_{L^2}^2d\tau\leq \frac12 \left(\|u_0\|_{L^2}^2+\|\te_0\|_{L^2}^2\right)e^t.
\end{equation}
Now, we apply the operator $\les$ to the velocity equation to find
\begin{equation*}
\begin{aligned}
\frac12\dt\|\les u\|_{L^2}^2+\|\Lambda^{1+\gamma+s} u\|_{L^2}^2&=-\int \Lambda^{1+\gamma+s}u\cdot\Lambda^{-1+\gamma+s}(u\cdot \grad u)+\int \Lambda^{2(\gamma+s)}u_2\te\\
&\leq\|\Lambda^{1+\gamma+s}u\|_{L^2}\|\Lambda^{-1+\gamma+s}(u\cdot \grad u)\|_{L^2}+\|\Lambda^{2(\gamma+s)}u\|_{L^2}\|\te\|_{L^2}.
\end{aligned}
\end{equation*}
Adding the energy inequality \eqref{energy} to the above one and noting that $2(s+\gamma)<s+\gamma+1$ for $s+\gamma\in(0,1)$ leads to
\begin{equation*}
\begin{split}
\frac12 \dt \|u\|_{H^{\gamma+s}}^2+\|\grad u\|_{L^2}^2+\| u\|_{H^{1+\gamma+s}}^2&\leq \|u\|_{L^2}^2+\|\te_0\|_{L^2}\|u\|_{H^{1+\gamma+s}}+\|\Lambda^{1+\gamma+s}u\|_{L^2}\|\Lambda^{-1+ \gamma+s}(u\cdot\grad u)\|_{L^2}.
\end{split}
\end{equation*}
By Sobolev embeddings and H\"older's inequality the last term above can be bounded as
\begin{equation*}
\|\Lambda^{-1+\gamma+s}(u\cdot\grad u)\|_{L^2}\leq c\|u\cdot \grad u\|_{L^{2/(2-\gamma-s)}}\leq c\|u\|_{L^{2/(1-\gamma-s)}}\|\grad u\|_{L^2}\leq c\|\les u\|_{L^2}\|\grad u\|_{L^2},
\end{equation*}
so by Young's inequality and \eqref{uest} we get
\begin{equation}\label{balance}
\frac12\dt\|u\|_{H^{\gamma+s}}^2+\frac12\|u\|_{H^{1+\gamma+s}}^2\leq \left(\|\te_0\|_{L^2}^2+\|u_0\|_{L^2}^2\right)e^t+c\|\grad u\|_{L^2}^2\| u\|_{H^{\gamma+s}}^2.
\end{equation}
To conclude, from \eqref{graduest} and Gr\"onwall's inequality applied to \eqref{balance}
it follows that
\begin{equation*}\label{usest}
\|u(t)\|_{H^{\gamma+s}}^2\leq c_1(\|u_0\|_{L^2},\|u_0\|_{H^{\gamma+s}},\|\te_0\|_{L^2},t),
\end{equation*}
\begin{equation*}\label{sgradusest}
\int_0^t\|u(t)\|_{H^{1+\gamma+s}}^2\leq c_2(\|u_0\|_{L^2},\|u_0\|_{H^{\gamma+s}},\|\te_0\|_{L^2},t).
\end{equation*}
These estimates can be justified by the usual limiting procedure \cite{MABE}.
\vspace{0.3cm}

Finally, to show uniqueness in this class, we consider two different solutions $(u_1,\theta_1)$, $(u_2,\theta_2)$ with the same initial data, denote the difference $\tilde{u}=u_2-u_1$, $\tilde{\theta}=\theta_2-\theta_1$ and take inner product:
\begin{equation*}
\begin{aligned}
\frac12\frac{d}{dt}\|\tilde{u}\|_{L^2}^2&\leq \|\grad u_2\|_{L^{\infty}}\|\tilde{u}\|_{L^2}^2-\|\grad\tilde{u}\|_{L^2}^2+\int\tilde{u}_2\tilde{\theta}dx\\
&\leq \|\grad u_2\|_{L^{\infty}}\|\tilde{u}\|_{L^2}^2-\|\grad\tilde{u}\|_{L^2}^2-\int \nabla \tilde{u}_2\cdot \nabla\Delta^{-1}\tilde{\theta}dx\leq \|\grad u_2\|_{L^{\infty}}\|\tilde{u}\|_{L^2}^2+2 \| \nabla\Delta^{-1}\tilde{\theta}\|_{L^2}^2,
\end{aligned}
\end{equation*}
where the incompressibility condition of $u_1$ and $u_2$ and integration by parts have been repeatedly used. As $\|\grad\Delta^{-1}\theta\|_{L^2}=\|\Lambda^{-1}\theta\|_{L^2}$, proceeding as above
\begin{equation*}
\begin{aligned}
\frac12\frac{d}{dt}\|\Lambda^{-1}\tilde{\theta}\|_{L^2}^2&=\int (\Delta^{-1}\tilde{\theta}) \grad\cdot(u_1 \tilde{\theta})dx+\int(\Delta^{-1}\tilde{\theta})\grad\cdot(\tilde{u}\theta_2) dx \\
&=-\int \grad \Delta^{-1}\tilde{\theta}\cdot \grad u_1 \cdot\grad\Delta^{-1}\tilde{\theta}-\int \grad \Delta^{-1} \tilde{\theta}\cdot \tilde{u} \theta_2 dx\\
&\leq \|\grad\Delta^{-1}\tilde{\theta}\|_{L^2}^2\|\grad u_1\|_{L^\infty}+\frac12\|\theta_2\|_{L^\infty} \|\tilde{u}\|_{L^2}^2+\frac12\|\theta_2\|_{L^\infty}\|\grad\Delta^{-1}\tilde{\theta}\|_{L^2}^2.
\end{aligned}
\end{equation*}
From the sum of both inequalities and Gr\"onwall's inequality the uniqueness follows.

\vspace{0.2cm}
Now the regularity of the velocity will be improved by using the vorticity equation. We will show that $u\in L^1(0,T;H^2)\cap L^1(0,T;C^{1+\gamma+\tilde{s}})$, $0<\tilde{s}<s$. We consider the vorticity equation as a forced heat equation
\begin{equation*}
\begin{aligned}
\omega_t-\Delta \omega=-u\cdot\grad \omega+\dxu\te,\hspace{0.5cm}\omega|_{t=0}=\omega_0,
\end{aligned}
\end{equation*}
and split the solution into three parts
\begin{equation}\label{vorsplitdef}
\omega(t)= \omega_1+\omega_2+\omega_3,\hspace{0.5cm}\omega_1=e^{t\Delta}\omega_0,\hspace{0.5cm}\omega_2=-\heatinv \nabla\cdot(u\omega),\hspace{0.5cm}\omega_3=\heatinv\partial_1\te.
\end{equation}

Since $\omega_0\in H^{-1+\gamma+s}$, by standard estimates for solutions of the heat equation in Sobolev spaces, we deduce that $\omega_1\in L^1(0,T;H^{1+\gamma+\tilde{s}})$ for $0<\tilde{s}<s$, and therefore $\omega_1\in L^1(0,T;C^{\gamma+\tilde{s}})$ by Sobolev embedding. This is sharp as in general $\omega_1\notin L^1(0,T;H^{1+\gamma+s})$ \cite{FERO2}. 

For the second part it suffices to prove that $u\omega\in L^q(0,T;H^{\gamma+s})$ for some $q>1$. In that case, one can make use of maximal regularity results for the heat equation to conclude that $\omega_2\in L^q(0,T;H^{1+\gamma+s})$ and hence $\omega_2\in L^1(0,T;C^{\gamma+s})$ (see for example \cite{FERO2}). From the fact that  $u\in L^\infty(0,T;H^{\gamma+s})\cap L^2(0,T;H^{1+\gamma+s})$, by interpolation it follows that, 
$$u\in L^\rho(0,T;H^{1+r}),\hspace{0.3cm}\rho=\frac{2}{1-(s+\gamma)+r}>2,\hspace{0.3cm}0<r<s+\gamma.$$
Since $\omega\in L^2(0,T;H^{\gamma+s})$ and $\|u\omega\|_{H^{\gamma+s}}\leq c\|u\|_{H^{1+r}}\|\omega\|_{H^{\gamma+s}}$, by H\"older's inequality it follows that 
$$u\omega\in L^q(0,T;H^{\gamma+s}),\hspace{0.2cm} \frac{1}{q}=\frac{1}{\rho}+\frac12 <1.$$

The temperature term, $\omega_3$, must be treated in a different way, as up to now $\te(t)$ only belongs to $L^2(\R^2)$, so the best we can achieve in the scale of Sobolev spaces is $\omega_3\in L^\infty(0,T;H^1)$. Anyway, even if we assume that $\te$ remains as a patch, it would belong to $H^\sigma$, $\sigma<1/2$ (as the characteristic function of a set with regular boundary belongs to these Sobolev spaces (see e.g. \cite{MAZYA})), but $\omega_3\in L^\infty(0,T;H^{1+\sigma})$ is not enough if $\gamma \geq 1/2$. We use instead that $\te\in L^\infty(0,T;L^\infty)\hookrightarrow  L^\infty(0,T;B_{\infty,\infty}^0)$. Then  $\omega_3\in L^\infty(0,T;B_{\infty,\infty}^1)$.
Since $C^r\hookrightarrow B_{\infty, \infty}^1$ $\forall r<1$, summing up the three terms and recalling that the Riesz transforms are continuous on Sobolev and H\"older spaces, we finally get the result $u\in L^1(0,T;H^2)\cap L^1(0,T;C^{1+\gamma+\tilde{s}})$.

\vspace{0.2cm}
\begin{rem}\label{remarktemp}
By interpolation, we could improve the time integrability to get $u\in L^p(0,T;H^2)$ for some $p=p(\gamma,s,\tilde{s})\in(1,2)$, although we won't use it.
However, to improve the spatial regularity, the restriction comes from the temperature term. Despite $\te_0\in H^{\sigma}$ for any $\sigma<1/2$, results in \cite{KUKA} are not sufficient to ensure that $\te$ remains in $H^\sigma$, due to the fact that the initial velocity has a low regularity comparable to that of the initial temperature. Once we prove that the patch is preserved and remains regular, the lines above would immediately imply that $u\in L^1(0,T;H^{2+\mu})$, $\mu<\{1/2,\gamma+s\}$.
\end{rem}
\vspace{0.2cm}

To conclude, as $\grad u\in L^1(0,T;C^\gamma)$ it follows that
\begin{equation*}
\|\grad X\|_{C^\gamma}\leq  \|\grad X_0\|_{C^\gamma}e^{\int_0^t\|\grad u\|_{L^\infty}d\tau}+\int_0^t\|\grad u(\tau)\|_{C^\gamma}\|\grad X(\tau)\|_{L^\infty}^{1+\gamma}e^{\int_\tau^t\|\grad u\|_{L^\infty}ds},
\end{equation*}
and hence $\|z\|_{L^\infty(0,T;C^{1+\gamma})}\leq C(T).$
\qed

\section{Control of curvature}\label{mainthe}

This section deals with the main result of this paper: the control of the curvature of the patch. To that end we bound $\grad^2 u$ in $L^1(0,T;L^\infty)$, which requires proving that the operators $\partial_k\partial_1(\partial_t-\Delta)^{-1}_0$ and $\partial_i^\perp\partial_j(-\Delta)^{-1}\partial_k\partial_1(\partial_t-\Delta)^{-1}_0$ applied to a patch are bounded.

\vspace{0.3cm}
\begin{thm}\label{Case2}
Let $D_0\subset \R^2$ be a bounded simply connected domain with boundary $\partial D_0\in W^{2,\infty}$, $u_0\in H^{1+s}$ with $s\in(0,1/2)$ a divergence-free vector field and $\te_0=1_{D_0}$ the characteristic function of $D_0$. Then, there exists a unique global solution $(u,\te)$ of \eqref{temperature}-\eqref{Boussinesq} such that
$$\te(x,t)=1_{D(t)}(x) \hspace{0.2cm}{\rm{and}}\hspace{0.2cm} \partial D \in L^\infty(0,T;W^{2,\infty})$$
 where $D(t)=X(D_0,t)$.
Moreover,
$$
u\in C(\R_+;H^{1+s})\cap L^2(0,T;H^{2+s})\cap L^p(0,T; H^{2+\mu})\cap L^q(0,T;W^{2,\infty}),
$$
for any $T>0$, $\mu<\frac12$, $1\leq p<2/(1+\mu-s)$, $1\leq q<2/(2-s)$.  
\end{thm}

Proof: From the previous theorem $\te$ remains as a $C^{1+\gamma}$ patch. Since the characteristic function of a set with regular boundary belongs to $H^\mu$ for any $\mu<1/2$, the result in \cite{KUKA} yields the existence and uniqueness of solutions $u\in C(\R_+;H^{1+s})\cap L^2(0,T;H^{2+s})$, $\te\in C(\R_+; H^s)$. 
\par
Proceeding as in the second part of the previous theorem, we split the vorticity as in \eqref{vorsplitdef}.
First, as $u\in L^\infty(0,T;H^{1+s})$ and $\omega\in L^2(0,T;H^{1+s})$ it follows  that $u\omega\in L^2(0,T;H^{1+s})$. Using this, by the properties of the heat kernel the first two parts are treated as before:
\begin{equation}\label{datoinicial}
\omega_0\in H^s \Longrightarrow \omega_1 \in L^1(0,T;H^{2+\tilde{s}}), \tilde{s}\in (0,s),
\end{equation}
\begin{equation}\label{nolineal}
u\omega\in L^2(0,T;H^{1+s})\Longrightarrow \omega_2\in L^2(0,T;H^{2+s}).
\end{equation}
 As the third term remains the same, $\omega_3\in L^\infty(0,T;H^{1+\mu})$, one obtains by interpolation  $u\in L^p(0,T;H^{2+\mu})$ where $p=\frac{2(1-(s-\tilde{s}))}{1+\mu-s-(s-\tilde{s})}>\frac43$.

We are thus left to prove that $u\in L^1(0,T;W^{2,\infty})$, as the persistence of regularity for the patch will follow by applying Gr\"onwall's inequality to the equation satisfied by the second derivatives of the particle trajectories:
$$\dt \grad^2X =\grad X^T\cdot \grad^2u (X)\cdot \grad X + \grad^2X\cdot \grad u(X). $$

\vspace{0.2cm}
In order to prove that $\grad^2 u \in L^1(0,T;L^\infty)$, we use the Biot-Savart formula
\begin{equation*}
\left(\grad^2 u\right)_{ijk}= \partial_k\partial_j u_i =-\partial_k(-\Delta)^{-1}\partial_i^{\perp}\partial_j (\omega_1+\omega_2+\omega_3),
\end{equation*}
where $i,j,k\in\{1,2\}$ and $\partial_1^{\perp}=-\partial_2$, $\partial_2^{\perp}=\partial_1$. Denote
\begin{equation}\label{hesu}
\begin{aligned}
(\grad^2 v_1)_{ijk}&=-\partial_k(-\Delta)^{-1}\partial_i^{\perp}\partial_j \omega_1=-\partial_k(-\Delta)^{-1}\partial_i^{\perp}\partial_j e^{t\Delta}\omega_0,\\
(\grad^2 v_2)_{ijk}&=-\partial_k(-\Delta)^{-1}\partial_i^{\perp}\partial_j \omega_2=\partial_k(-\Delta)^{-1}\partial_i^{\perp}\partial_j \heatinv\grad\cdot( u\omega),\\
(\grad^2 v_3)_{ijk}&=-\partial_k(-\Delta)^{-1}\partial_i^{\perp}\partial_j \omega_3=-\partial_k(-\Delta)^{-1}\partial_i^{\perp}\partial_j \heatinv\partial_1\te.
\end{aligned}
\end{equation}
As $\omega_1 \in L^1(0,T;H^{2+\tilde{s}})$ and $\omega_2 \in L^2(0,T;H^{2+s})$, it is clear that $\grad^2 v_1,\grad^2 v_2 \in L^1(0,T;L^\infty)$.

\vspace{0.2cm}
The temperature term can be written as an operator applied to $\te$: 
\begin{equation}\label{nabla2u3}
\left(\grad^2 v_3\right)_{ijk}= S_{ijk} \te = R_iR_j\partial_k\partial_1 \heatinv \te
\end{equation}
where $R_i$, $R_j$ denote Riesz transforms.

\vspace{0.3cm}
The rest of the proof is structured as follows: first, in Section $3.1$ we defined the operators $\partial_k\partial_1\left(\partial_t-\Delta\right)^{-1}_0$ and $R_iR_j\partial_k\partial_1\left(\partial_t-\Delta\right)^{-1}_0$. Later in Section \ref{bounds} we proceed to bound these operators applied to the patch, i.e., we compute the bounds in $L^1(0,T;L^\infty)$ for the gradient of the vorticity and the second derivatives of the velocity. The bounds found give $\grad^2 u\in L^q(0,T; L^\infty)$, $1\leq q<2/(2-s)$ (see Remark \ref{remarkcurvature}).

\subsection{Definition of the operators $S_{ijk}$}

\vspace{0.4cm}
\subsubsection{Singular heat kernels:
	 $\partial_k\partial_1 \heatinv$}

Denote the heat kernel $K(x,t)=\frac{1}{4\pi t}e^{-|x|^2/4t}$. Then, the operators for the patch can be written as
\begin{equation}\label{identidad}
\partial_k \omega_3(x,t)=\partial_k\partial_1 \heatinv \te (x,t):=\lim_{\eps\to 0} \int_0^{t-\eps}\int_{\R^2}(\partial_k\partial_1 K)(x-y,t-\tau)\te(y,\tau)dy d\tau,
\end{equation}
where
\begin{equation}\label{nucleos}
\left\{\begin{aligned}
&\partial_1^2K(x,t)=\frac{1}{8\pi t^2}\left(\frac{x_1^2}{2t}-1\right)e^{-|x|^2/4t},\\
&\partial_2\partial_1 K(x,t)=\frac{x_1x_2}{16\pi t^3}e^{-|x|^2/4t}.
\end{aligned}\right.
\end{equation}

\medskip
Note that the principal value is needed: due to the singularity of the kernels at $(x,t)=(0,0)$  along parabolas $\tau=M r^2$ ($M\neq 1/4$ for $\partial_1^2K$), the kernels are not absolutely integrable:

\begin{equation*}
\int_0^t \int_{\R^2}|\partial_1^2K(y,\tau)|dy d\tau=+\infty.
\end{equation*}
It is  clear then that these operators are not bounded for a general $L^\infty$ function. 

Here we point out that while $\partial_1\partial_2K(\cdot,t)$ clearly has  zero mean on circles for any $t$, this is not true for $\partial_1^2K$. However, the time integral provides a new kind of cancellation:
\begin{equation}\label{cancelacion}
\begin{aligned}
\lim_{\eps\to 0}\int_\eps^t \int_{B_R}\partial_1^2K(y,\tau)dyd\tau&=\lim_{\eps\to 0}\int_\eps^t \int_0^R\frac{r}{4\tau^2}\left(\frac{r^2}{4\tau}-1\right)e^{-r^2/4\tau}d\tau dr\\
&=\lim_{\eps\to0}\int_\eps^t \frac{-R^2e^{-R^2/4\tau}}{8\tau^2}d\tau =-\frac12e^{-R^2/4t},
\end{aligned}
\end{equation}
where $B_R$ is the ball of radius $R$.

\vspace{0.6cm}
\subsubsection{Oseen-type kernels: $R_iR_j\partial_k\partial_1  \heatinv$}

We will see that  $\partial_k\omega_3=\partial_k\partial_1\heatinv\te$ is bounded if $\te(t)$ is a $C^{1+\gamma}$ patch for all $t$. To get the boundedness of $R_iR_j\partial_k\partial_1 \heatinv$ we need to treat the combined operator directly, which we will write as a convolution with an explicit kernel:
\begin{equation*}
\begin{aligned}
\left(\grad^2v_3\right)_{ijk} (x,t)&=\frac{1}{2\pi}\int_{\R^2}\frac{y_i^\perp}{|y|^2}\partial_j\partial_k\omega_3(x-y,t)dy=\frac{1}{2\pi}\int_{\R^2}\log{|y|}\partial_i^\perp\partial_j\partial_k\omega_3(x-y,t)dy\\
&=\partial_i^\perp\partial_j\partial_k\partial_1 \lim_{\eps\to0}\int_0^{t-\eps}\int_{\R^2} \left(\frac{1}{2\pi} \int_{\R^2}K(x-y-z,t-\tau)\log{|y|}dy\right)\theta(z,\tau)dzd\tau.
\end{aligned}
\end{equation*}
We notice now that the term in brackets can be seen as the solution of the Laplace equation in $\R^2$ with force $K(x,t)$, which can be computed explicitly:
\begin{equation*}
\Delta^{-1}K(x,t):=\frac{1}{2\pi}\int\frac{1}{4\pi t}e^{-\frac{|x-y|^2}{4t}}\log|y| dy=\frac{1}{2\pi}\left(\log|x|+\frac12\int_{|x|^2/4t}^\infty \frac{e^{-z}}{z}dz\right).
\end{equation*}

\medskip
Thus the operators can be written as
\begin{equation*}
\begin{aligned}
R_iR_j\partial_k\partial_1 \heatinv\te(x,t)&=\partial_i^\perp\partial_j\partial_k\partial_1 \int_0^{t} (\Delta^{-1}K)(t-\tau)*\te(\tau)(x)d\tau,
\end{aligned}
\end{equation*}
or, as a kernel convolution in $\theta$
\begin{equation*}
\begin{aligned}
R_iR_j\partial_k\partial_1 \heatinv\te(x,t)&:=\lim_{\eps\to0} \int_0^{t-\eps} K_{ijk}(t-\tau)*\te(\tau)(x)d\tau,
\end{aligned}
\end{equation*}
where
\begin{equation}\label{kernelsref}
\begin{aligned}
K_{ijk}(x,t):=\partial_1\partial_j\partial_i^\perp\partial_k (\Delta^{-1}K(x,t))&=\partial_1\partial_j\partial_i^\perp\partial_k\left(\frac{1}{2\pi}\int\frac{1}{4\pi t}e^{-\frac{|x-y|^2}{4t}}\log|y| dy\right)\\
&=\frac{1}{2\pi}\partial_j\partial_i^\perp\partial_k \left(\frac{x_1}{|x|^2}\left(1-e^{-|x|^2/4t}\right)\right).
\end{aligned}
\end{equation}

\medskip
The eight possible kernels reduce to four. Here we show the expression of $K_{112}$ as an example (see the Appendix for the rest):
\begin{equation}\label{kernels}
\begin{aligned}
K_{112}&=\left(\frac{24x_1^2x_2^2}{\pi |x|^5}-\frac{3}{\pi|x|}\right)G(|x|,t)-\frac{e^{-|x|^2/4t}}{\pi(4t)^2}\left(-2+\frac{4}{4t}\frac{x_1^2x_2^2}{|x|^2}+12\frac{x_1^2x_2^2}{|x|^4}\right),
\end{aligned}
\end{equation}
where
\begin{equation}\label{Gfunction}
G(|x|,t)=\frac{1}{|x|^3}(1-e^{-|x|^2/4t})-\frac{1}{4t |x|}e^{-|x|^2/4t}.
\end{equation}

All these kernels are not integrable, but they show  cancellations in circles. In fact, the kernels $K_{111}$ and $K_{122}$ have zero mean on circles, while for $K_{112}$ and $K_{211}$ we find the cancellation \eqref{cancelacion}:
\begin{equation}\label{cancelacionK112}
\begin{aligned}
\lim_{\eps\to0}\int_\eps^{t}\int_{B_R} K_{112}(y,\tau)dyd\tau=\lim_{\eps\to0}\int_\eps^{t}\int_0^R\frac{r}{(4\tau)^2}\left(\frac{r^2}{4\tau}-1\right)e^{-r^2/4\tau}drd\tau=-\frac12e^{-R^2/4t},\\
\lim_{\eps\to0}\int_\eps^{t}\int_{B_R} K_{211}(y,\tau)dyd\tau=\lim_{\eps\to0}\int_\eps^{t} \int_0^R \frac{3r}{(4\tau)^2}\left(1-\frac{r^2}{4\tau}\right)e^{-r^2/4\tau}drd\tau=-\frac32e^{-R^2/4t}.
\end{aligned}
\end{equation}
We split these kernels in two parts,  $K_{ijk}=K_{ijk}^*+K_{ijk}^o$, where $\int_{\partial B_R}K_{ijk}^o(y)d\sigma(y)=0$.

\subsection{Bound for $\nabla^2 u_3$}\label{bounds}

In this section we prove that the gradient of the vorticity and the second derivatives of the velocity are bounded in $L^1(0,T;L^\infty)$ if $\theta$ is a patch.
First we prove that $\partial_k\omega_3\in L^1(0,T;L^\infty)$.
We do it carefully for $\partial_1\omega_3$ and then show the differences with $\partial_2\omega_3$.  

Let $\gamma\in(0,1)$ and consider a parametrization $z(\cdot,t)\in C^{1,\gamma}$ as in \eqref{parametrization}, which we know that exists thanks to Theorem \ref{Case1}  and verifies
$$|\partial_\alpha z|_{{\rm{inf}}}(t)=\inf_{\alpha\in [0,1]}|\partial_\alpha z(\alpha,t)|>0,$$
$$|\partial_\alpha z|_\gamma(t)=\sup_{\alpha\neq \beta}\frac{|\partial_\alpha z(\alpha,t)-\partial_\alpha z(\beta,t)|}{|\alpha-\beta|^\gamma}<\infty,$$
for all $t\in[0,T]$.
Take the cut-off distance
\begin{equation*}\label{deltadef}
\delta=\min_{\tau\in[0,t]}\left(\frac{|\partial_\alpha z|_{{\rm{inf}}}(t)}{|\partial_\alpha z|_\gamma(t)}\right)^{1/\gamma}>0,
\end{equation*}
which is a fixed positive quantity by the previous result.
Denote 
$$d(x,t):=d(x,\partial D(t))$$
 the distance from $x$ to $\partial D(t)$ and write
\begin{equation}\label{domega}
\begin{gathered}
\partial_1\omega_3(x,t)= \pv \int_0^t \int_{\R^2} \frac{1}{8\pi (t-\tau)^2}\left(\frac{(x_1-y_1)^2}{2(t-\tau)}-1\right)e^{-|x-y|^2/4(t-\tau)}\te(y,\tau) dy d\tau=I_1+I_2,
\end{gathered}
\end{equation}
where $\textup{pv}$ denotes principal value defined as in \eqref{identidad} and the splitting consists in
\begin{equation*}
I_1=  \int_0^t \int_{D(t-\tau)\cap |y|\geq \delta} \frac{1}{8\pi \tau^2}\left(\frac{y_1^2}{2\tau}-1\right)e^{-|y|^2/4\tau} dy d\tau,
\end{equation*}
\begin{equation*}\label{I2}
I_2=  \pv \int_0^t \int_{D(t-\tau)\cap |y| \leq \delta} \frac{1}{8\pi \tau^2}\left(\frac{y_1^2}{2\tau}-1\right)e^{-|y|^2/4\tau} dy d\tau.
\end{equation*}
The bound for $I_1$ follows from the fact that $\delta$ is a fixed positive quantity:
\begin{equation}\label{boundI1}
\begin{aligned}
|I_1|&\leq \int_0^t \int_\delta^\infty \int_0^{2\pi} \left(\frac{r^3\cos^2\alpha}{2\tau^3}+\frac{r}{\tau^2}\right)\frac{e^{-r^2/4\tau}}{8\pi} d\alpha dr d\tau=\left(\frac{4t}{\delta^2}+\frac12\right)e^{-\delta^2/4t}.
\end{aligned}
\end{equation}
Now we deal with $I_2$. First, if $\delta< d(x,\tau)$ for all $\tau\in[0,t]$, then cancellation \eqref{cancelacion} yields
\begin{equation}\label{boundI2facil}
|I_2|\leq\frac12 e^{-\delta^2/4t}.
\end{equation}
Let's consider then the case in which the boundary of the patch is close to $x$, $0<d(x,\tau)\leq \delta$. Since the velocity of the patch is finite almost everywhere in time and space ($u\in L^\infty(0,T;L^\infty)$), then distance $d(x,\tau)$ is a Lipschitz function in time in $[0,T]$. In fact,

\begin{equation*}
\begin{aligned}
d(x,\tau)=\|x-x^*(\tau)\|,\hspace{0.2cm} x^*(\tau)=\underset{y(\tau)\in\partial D(\tau)}{{\rm{arg \hspace{0.1cm}min}}}\|x-y(\tau)\|,
\end{aligned}
\end{equation*}
\begin{equation*}\label{intdistancia}
\begin{aligned}
 d(x,\tau)=d(x,0)+\int_0^\tau u(x^*(\tau'))\cdot n(x^*(\tau'))d\tau',
\end{aligned}
\end{equation*}
so
\begin{equation*}\label{vel}
|d(x,\tau)-d(x,0)|\leq U \tau, \hspace{0.4cm} |d(x,t)-d(x,\tau)|\leq U(t-\tau),
\end{equation*}
where $U=\|u\|_{L^\infty(0,T;L^\infty)}$. If we denote
\begin{equation*}
0<\eps= \min_{\tau\in[0,t]}d(x,\tau)\leq d(x,t)\leq \delta,
\end{equation*}
we may consider two possible cases:
\vspace{0.2cm}
 \begin{equation}\label{cases12}
 \textup{\textbf{\em{Case 1}:}}
  \lim_{\eps\to 0^+}\sup_{\eta\in[\eps,\delta]}\frac{d(x,t)}{\eta}\leq 2,\hspace{1cm}\textup{\textbf{\em{Case 2}}}:  \lim_{\eps\to 0^+}\sup_{\eta\in[\eps,\delta]}\frac{d(x,t)}{\eta}\geq 2.
 \end{equation}
  \vspace{0.2cm}
\medskip
\newline
In both cases the evolution of the distance satisfies that
\begin{equation}\label{distbound}
\epsilon\leq d(x,t-\tau)\leq d(x,t)+U\tau \hspace{0.6cm}\forall \tau \in [0,t].
\end{equation}
In the second case, we will also use that for $\tau\in[0,(d(x,t)-\epsilon)/2U]$, 
\begin{equation}\label{case22}
\frac{d(x,t)}{2} \leq d(x,t-\tau)\leq \frac32 d(x,t).
\end{equation}
\medskip
In all cases we decompose the space domain $D(\tau)\cap |x-y|\leq \delta$ as follows:
\begin{equation}\label{descomposicion}
\begin{gathered}
S_r(x,\tau)=\{z\in \R^2: |z|=1, x+rz\in D(\tau)\},\\
\Sigma(x,\tau)=\{z\in\R^2:|z|=1, \grad \varphi(x^*(\tau))\cdot z\geq 0)\},\\
R_r(x,\tau)=(S_r(x,\tau)\backslash \Sigma(x,\tau))\cup (\Sigma(x,\tau)\backslash S_r(x,\tau) ),
\end{gathered}
\end{equation}
and use for each $\tau$ the {\em{Geometric Lemma}}  in \cite{CONST}
\begin{equation}\label{geometriclema}
|R_r(x,\tau)|\leq 2\pi \left((1+2^\gamma)\frac{d(x,\tau)}{r}+2^\gamma\frac{r^\gamma}{\delta^\gamma}\right),
\end{equation}
valid for all $\tau\in[0,T]\medskip$. From \eqref{I2} we get
\begin{equation}\label{splitJ}
\begin{aligned}
|I_2|&\leq J_1+J_2+J_3,\hspace{0.4cm}\textup{where}\hspace{0.4cm}
J_1=\left|\pv\int_0^t \int_0^{d(x,t-\tau)}\frac{r}{4\tau^2}\left(\frac{r^2}{4\tau}-1\right)e^{-r^2/4\tau}drd\tau \right|,\\
&J_2=\left| \int_0^t \int_{d(x,t-\tau)}^\delta \int_{\Sigma(x,t-\tau)} \frac{r}{8\pi\tau^2}\left(\frac{r^2\cos^2\alpha}{2\tau}-1\right)e^{-r^2/4\tau} d\alpha dr d\tau\right|,\\
&J_3=\left| \int_0^t \int_{d(x,t-\tau)}^\delta \int_{R_r(x,t-\tau)} \frac{r}{8\pi\tau^2}\left(\frac{r^2\cos^2\alpha}{2\tau}-1\right)e^{-r^2/4\tau} d\alpha dr d\tau\right|.
\end{aligned}
\end{equation}

\medskip
\underline{{\em{Case 1}}}: In $J_1$ the time-space cancellation in  \eqref{cancelacion} and \eqref{distbound} yields
\begin{equation}\label{boundJ1}
\begin{aligned}
J_1&\leq \int_0^t \frac{d(x,t-\tau)^2}{8\tau^2}e^{-d(x,t-\tau)^2/4\tau}d\tau\leq
\int_0^t \left(\frac{d(x,t)^2}{4\tau^2}+\frac{U^2}{4}\right)e^{-\eps^2/4\tau}d\tau\\
&\leq \frac{d(x,t)^2}{\eps^2}e^{-\eps^2/4t}+\frac{U^2}{4}t\leq 4+\frac{U^2}{4}t.
\end{aligned}
\end{equation}
By the parity of the kernel, $J_2$ can be estimated in a similar way to $J_1$:
\begin{equation}\label{boundJ2}
\begin{aligned}
J_2&=\left|\frac12\int_0^t\left( \frac{-\delta^2}{8\tau^2}e^{-\delta^2/4\tau}+\frac{d(x,t-\tau)^2}{8\tau^2}e^{-d(x,t-\tau)^2/4\tau}\right)d\tau  \right|\leq 2+\frac{U^2}{8}t.
\end{aligned}
\end{equation}
Finally, to bound $J_3$ we use the geometric lemma \eqref{geometriclema} and again \eqref{distbound}
\begin{equation}\label{J3}
\begin{aligned}
J_3&\leq \int_0^t \int_{d(x,t-\tau)}^\delta \frac{r}{8\pi\tau^2}\left(\frac{r^2}{2\tau}+1\right)|R_r(x,t-\tau)|e^{-r^2/4\tau}drd\tau\\
&\leq 3\int_0^t\int_{d(x,t-\tau)}^\delta\frac{r}{8\pi\tau^2}|R_r(x,t-\tau)|e^{-r^2/8\tau}drd\tau\leq L_1+L_2,
\end{aligned}
\end{equation}
where
\begin{equation}\label{boundL2}
\begin{aligned}
L_1&=3\frac{1+2^\gamma}{4}\int_0^t \int_{\eps}^\delta \frac{(d(x,t)+U\tau)}{\tau^2}e^{-r^2/8\tau}drd\tau\leq 9d(x,t)\int_{\eps}^\delta\frac{2}{r^2}e^{-r^2/8t}dr\\
&\quad+9U\int_0^t \frac{\sqrt{2}}{2\sqrt{\tau}}d\tau\int_0^\infty e^{-y^2}dy
\leq 18\frac{d(x,t)}{\eps}+9\frac{U\sqrt{2\pi}}{2}\sqrt{t}\leq 36+\frac{9U\sqrt{2\pi}}{2}\sqrt{t},\\
L_2&= 3\frac{2^\gamma}{\delta^\gamma}\int_0^t \int_{\eps}^\delta \frac{r^{1+\gamma}}{4\tau^2}e^{-r^2/8\tau}drd\tau\leq 3\frac{2^\gamma}{\delta^\gamma}\int_{\eps}^\delta\frac{2}{r^{1-\gamma}}dr\leq 6\frac{2^\gamma}{\gamma}.
\end{aligned}
\end{equation}
Joining the above bounds we finally find that
\begin{equation}\label{boundJ3}
J_3\leq 36+6\frac{2^\gamma}{\gamma}+\frac{9U\sqrt{2\pi}}{2}t^{1/2}.
\end{equation}
From the splitting \eqref{splitJ} and the bounds \eqref{boundJ1}, \eqref{boundJ2} and \eqref{boundJ3}, it is easy to get
\begin{equation}\label{boundI2}
|I_2|\leq 42+6\frac{2^\gamma}{\gamma}+\frac{9U\sqrt{2\pi}}{2}t^{1/2}+\frac{3U^2}{8}t.
\end{equation}
\medskip
Thus, from \eqref{domega} and the bounds \eqref{boundI1}, \eqref{boundI2facil} and \eqref{boundI2}, we conclude that
\begin{equation*}
|\partial_1\omega_3(x,t)|\leq \frac{85}{2}+6\frac{2^\gamma}{\gamma}+\frac{9U\sqrt{2\pi}}{2}t^{1/2}+\left(\frac{3U^2}{8}+\frac{4}{\delta^2}\right)t,
\end{equation*}
and therefore 
\begin{equation}\label{boundhola}
\begin{aligned}
\|\partial_1\omega_3\|_{L_T^1(L^\infty)}\leq \left(\frac{85}{2}+6\frac{2^\gamma}{\gamma}\right)T+\frac{3U\sqrt{2\pi}}{2}T^{3/2}+\left(\frac{3U^2}{16} +\frac{2}{\delta^2}\right)T^2.
\end{aligned}
\end{equation}

\vspace{0.3cm}
\medskip
\underline{{\em{Case 2}}}: To get the bounds now we need to split the time integrals in $[0,t^*]$ and $[t^*,t]$, with $t^*=\frac{d(x,t)-\epsilon}{2U}>0$ and use the Lipschitz continuity of the velocity to remove the singularity at time $t$. In this way, for the first interval one has the bounds $\frac12d(x,t)\leq d(x,t-\tau)\leq \frac32 d(x,t)$. Proceeding as in \eqref{splitJ}, we now split further in two time intervals and use \eqref{case22} for the first and \eqref{distbound} for the second:
\begin{equation*}\label{splitJ1case2}
\begin{aligned}
J_1&\leq \int_0^t \frac{d(x,t-\tau)^2}{8\tau^2}e^{-d(x,t-\tau)^2/4\tau}d\tau\\
&\leq \frac94 \int_0^{t^*}\frac{d(x,t)^2}{8\tau^2}e^{-d(x,t)^2/16\tau}d\tau+\int_{t^*}^t\frac{(d(x,t)+U\tau)^2}{8\tau^2}e^{-\eps^2/4\tau}d\tau=L_3+L_4.
\end{aligned}
\end{equation*}
As $t^*=(d(x,t)-\eps)/2U>0$, one easily gets 
\begin{equation*}\label{boundJ11}
\begin{aligned}
L_3&\leq \frac92 e^{-d(x,t)^2/16t^*}\leq \frac92.
\end{aligned}
\end{equation*}
For $L_4$ we proceed as follows
\begin{equation}\label{boundJ12}
\begin{aligned}
L_4&\leq\int_{t^*}^t \frac{d(x,t)^2}{4\tau^2}e^{-\epsilon^2/4\tau}d\tau +\int_{t^*}^t\frac{U^2}{4}e^{-\epsilon^2/4\tau}d\tau\leq\left(e^{-\eps^2/4t}-e^\frac{-\eps^2U/2}{d(x,t)-\eps}\right)\frac{d(x,t)^2}{\eps^2}+\frac{U^2}{4}t\\
&\leq \left(\frac{U/2}{d(x,t)-\eps}-\frac{1}{4t}\right)d(x,t)^2+\frac{U^2}{4}t\leq  \frac{U/2}{1-\eps/d(x,t)}\delta+\frac{U^2}{4}t\leq U\delta+\frac{U^2}{4}t,
\end{aligned}
\end{equation}
where we have used that in this case $\eps/d(x,t)\leq 1/2$.
 Thus the final bound for $J_1$ is given by
\begin{equation}\label{boundJ1case2}
J_1\leq \frac92+U\delta+\frac{U^2}{4}t.
\end{equation}
The second term can be treated similarly:
\begin{equation}\label{boundJ2case2}
J_2\leq\left| -\frac14e^{-\delta^2/4t}+\frac12\int_0^t\frac{d(x,t-\tau)^2}{8\tau^2}e^{-d(x,t-\tau)^2/4\tau}\right|d\tau\leq \frac52+\frac12U\delta+\frac{U^2}{8}t.
\end{equation}
Finally, we split the time integral in the $J_3$ term \eqref{splitJ}
\begin{equation*}
J_3\leq L_5+L_6,
\end{equation*}
where
\begin{equation*}
L_5=3\int_0^{t^*}\int_{d(x,t-\tau)}^\delta \frac{r}{4\tau^2}\left((1+2^\gamma)\frac{d(x,t-\tau)}{r}+\frac{2^\gamma}{\delta^\gamma}r^\gamma\right)e^{-r^2/8\tau}drd\tau,
\end{equation*}
\begin{equation*}
L_6=9\int_{t^*}^t\int_{d(x,t-\tau)}^\delta \frac{d(x,t-\tau)}{4\tau^2}e^{-r^2/8\tau}drd\tau+3\frac{2^\gamma}{\gamma}\int_{t^*}^t\int_{d(x,t-\tau)}^\delta\frac{r^{1+\gamma}}{4\tau^2}e^{-r^2/8\tau}drd\tau.
\end{equation*}
The term $L_5$ is done as $J_3$ in {\em{case 1}} \eqref{J3} but using the lower bound in \eqref{case22}. This yields
\begin{equation*}
L_5\leq \frac{27}{2}+6\frac{2^\gamma}{\gamma}.
\end{equation*}
The second term in $L_6$ is bounded by $L_2$ \eqref{boundL2}, while for the first we proceed as in \eqref{boundJ12}
\begin{equation*}
L_6\leq 9d(x,t)\int_\eps^\delta \frac{2}{r^2}\left(\frac{r^2}{8t^*}-\frac{r^2}{8t}\right)dr+9U\int_{t^*}^t \frac{\sqrt{2}}{2\sqrt{\tau}}d\tau\int_0^\infty e^{-y^2}dy\leq 9U\delta+\frac{9U\sqrt{2\pi}}{2}t^{1/2}.
\end{equation*}
Summing up the last two bounds above we obtain
\begin{equation}\label{boundJ3case2}
J_3\leq \frac{27}{2}+12\frac{2^\gamma}{\gamma}+9U\delta+\frac{9U\sqrt{2\pi}}{2}t^{1/2}.
\end{equation}
Finally, from \eqref{domega} and joining the above bounds \eqref{boundI1}, \eqref{boundI2facil}, \eqref{boundJ1case2}, \eqref{boundJ2case2} and \eqref{boundJ3case2}, we get 
\begin{equation*}
\begin{aligned}
|\partial_1\omega_3|&\leq 21+\frac{21}{2}U\delta+12\frac{2^\gamma}{2}+\frac{9U\sqrt{2\pi}}{2}t^{1/2}+\left(\frac{3U^2}{8}+\frac{4}{\delta^2}\right)t,
\end{aligned}
\end{equation*}
so after integration in time,
\begin{equation*}
\|\partial_1\omega_3\|_{L^1_T(L^\infty)}\leq \left(\frac{41}{2}+\frac{21}{2}U\delta+12\frac{2^\gamma}{2}\right)T+3U\sqrt{2\pi}T^{3/2}+\left(\frac{3U^2}{16}+\frac{2}{\delta^2}\right)T^2.
\end{equation*}
From the bound above and \eqref{boundhola}, the final bound for $\partial_1\omega_3$ in all cases is
\begin{equation}\label{bounddomega3}
\|\partial_1\omega_3\|_{L^1_T(L^\infty)}\leq c_1(\gamma, U,\delta)T+c_2(U)T^{3/2}+c_3(\delta,U)T^2.
\end{equation}

\vspace{0.3cm}
Note that for $\partial_2\omega_3$ the proof above reduces to bound the term $J_3$ in \eqref{splitJ}, as it has zero mean on half circles. If we write the corresponding kernels \eqref{nucleos} in polar coordinates
\begin{equation*}
|\partial_1^2K(r,\alpha,t)|\leq\frac{r^2}{8\pi t^3}e^{-r^2/4t}+\frac{1}{8\pi t^2}e^{-r^2/4t},\hspace{1cm}|\partial_1\partial_2K(r,\alpha,t)|\leq \frac{r^2}{16\pi t^3}e^{-r^2/4t},
\end{equation*}
one can see that the bound \eqref{bounddomega3} is also valid for $\partial_2\omega_3$.
\qed

\vspace{0.5cm}
We proceed now to bound the second derivatives of the velocity. We want to show that $R_iR_j\partial_k\omega_3\in L^1(0,T;L^\infty)$. The operators $K_{ijk}$ involved can be decomposed in two parts, one with zero mean on half circles and another with the cancellation \eqref{cancelacion}, as shown in \eqref{cancelacionK112}. Rewrite for example the kernel $K_{112}=K_{112}^o+K_{112}^*$ in \eqref{kernels} as
\begin{equation*}
K_{112}^*(r,\alpha)=\frac{1}{\pi(4t)^2}e^{-r^2/4t}\left(2-\frac{4r^2\cos^2\alpha\sin^2\alpha}{4t} -12\cos^2\alpha\sin^2\alpha\right),
\end{equation*}
\begin{equation*}
K_{112}^o(r,\alpha)=\frac{1}{\pi r}\left(3-24\cos^2\alpha \sin^2\alpha\right)G(r,t),
\end{equation*}
where $G(r,t)$ is given by \eqref{Gfunction}.

Using the same decomposition \eqref{descomposicion}, we notice from \eqref{cancelacionK112} that the part corresponding to $K_{112}^*$ can be estimated in the same way as we did with $\partial_1^2K$ obtaining the same bound (up to a constant). Following the splittings \eqref{domega} and \eqref{splitJ}, as $K_{112}^o$ has zero mean on half circles, we only need to estimate its contribution due to $R_r$:

\begin{equation*}
\begin{aligned}
J_4&=\left|\int_0^t\int_{d(x,t-\tau)}^\delta \int_{R_r(x,t-\tau)}r K_{112}^o(r,\alpha) d\alpha dr d\tau\right| \\
&\leq \int_0^t\int_{d(x,t-\tau)}^\delta\int_0^{2\pi}\frac{1}{\pi}|3-24\cos^2{\alpha}\sin^2{\alpha}|\left|G(r,\tau)\right|\left|R_r(x,t-\tau)\right|d\alpha drd\tau.
\end{aligned}\medskip
\end{equation*}
Using \eqref{geometriclema} and the fact that $G(r,\tau)\geq 0$ $\forall r,\tau\geq 0$, one gets
\begin{equation}\label{boundIgen}
J_4\leq 54 \int_0^t\int_{d(x,t-\tau)}^\delta G(r,\tau)\left((1+2^\gamma)\frac{d(x,t-\tau)}{r}+2^\gamma \frac{r^\gamma}{\delta^\gamma}\right)drd\tau\leq 54(L_7+L_8),
\end{equation}
where
\begin{equation*}
\begin{aligned}
L_7=3 \int_0^t\int_{d(x,t-\tau)}^\delta G(r,\tau)\frac{d(x,t-\tau)}{r}drd\tau,\hspace{0.7cm}
L_8=\frac{2^\gamma}{\delta^\gamma}\int_0^t\int_{d(x,t-\tau)}^\delta G(r,\tau) r^\gamma drd\tau.
\end{aligned}
\end{equation*}
In the term $L_8$ the singularity has been removed so we can integrate it directly
\begin{equation}\label{boundI2vel}
L_8=\frac{2^\gamma}{\delta^\gamma}\int_{0}^\delta t \frac{1-e^{-r^2/4t}}{r^{3-\gamma}}  dr \leq \frac{2^\gamma}{4\gamma}.
\end{equation}
To deal with the term $L_7$ we consider the {\em{cases 1}} and {\em{2}} \eqref{cases12} separately.
In {\em{case 1}}, we use first \eqref{distbound}  
\begin{equation*}
L_7\leq 3d(x,t)\int_0^t \int_\eps^\delta \frac{G(r,\tau)}{r}drd\tau+ 3U\int_0^t \int_\eps^\delta G(r,\tau)\frac{\tau}{r}drd\tau=M_1+M_2.
\end{equation*}
Again, in $M_2$ we have enough cancellation to integrate
\begin{equation}\label{boundM2}
M_2=\frac{U}{4}\int_{0}^t \left(\frac{1}{\sqrt{\tau}}\int_0^{\delta/2\sqrt{\tau}}e^{-y^2}dy+\frac{1}{\delta}e^{-\delta^2/4\tau}-4\tau\frac{1-e^{-\delta^2/4\tau}}{\delta^3}\right)d\tau\leq \frac{U\sqrt{\pi}}{4}t^{1/2}+\frac{U}{4\delta}t.
\end{equation}
Because we are dealing with {\em{case 1}} we find that 
\begin{equation*}
M_1\leq 3d(x,t)\int_\eps^\delta t\frac{1-e^{-r^2/4t}}{r^4}dr\leq \frac{3d(x,t)}{4\eps}\leq \frac32.
\end{equation*}
Joining the above bounds we find that in {\em{case 1}} 
\begin{equation*}
J_4\leq 81+\frac{2^\gamma}{4\gamma}+\frac{27U\sqrt{\pi}}{2}t^{1/2}+\frac{U}{2\delta}t.
\end{equation*}
In {\em{case 2}} we split the time integral $L_7=M_3+M_4$ where
\begin{equation*}
M_3=3\int_0^{t^*}\int_{d(x,t-\tau)}^\delta G(r,\tau)\frac{d(x,t-\tau)}{r}drd\tau,\hspace{0.4cm}M_4=3 \int_{t^*}^t\int_{d(x,t-\tau)}^\delta G(r,\tau)\frac{d(x,t-\tau)}{r}drd\tau.
\end{equation*}
Proceeding as in $M_1$ in {\em{case 1}}  but using the lower and upper bounds in \eqref{case22} we find that
\begin{equation*}
M_3\leq \frac92.
\end{equation*}
Finally, we split $M_4=N_1+N_2$, where
\begin{equation*}
N_1=3d(x,t)\int_{t^*}^t\int_\eps^\delta \frac{G(r,\tau)}{r}drd\tau,\hspace{0.6cm}N_2=3U\int_{t^*}^{t}\int_\eps^\delta G(r,\tau)\frac{\tau}{r}drd\tau.
\end{equation*}
The second term is bounded by $M_2$ in {\em{case 1}} \eqref{boundM2}. To deal with $N_1$ we need to take advantage of $t^*$ by applying the mean value theorem
\begin{equation*}
\begin{aligned}
N_1=3d(x,t)\int_\eps^\delta \frac{1}{r^4}\left(t\left(1-e^{-r^2/4t}\right) -t^*\left(1-e^{-r^2/4t^*}\right) \right)dr\leq 3d(x,t)\frac{\delta}{2}\frac{1}{4t^*}\leq \frac32 \delta U.
\end{aligned}
\end{equation*}
where we have use that in this case $\eps/d(x,t)\leq 1/2$. Therefore we conclude that in all cases the following bound holds:
\begin{equation}\label{boundmeanzero}
J_4\leq c_1(\delta, \gamma,U)+c_2(U)t^{1/2}+c_3(\delta,U)t.
\end{equation}
\newline
One can see from their expression in \eqref{allthekernels} that all the terms appearing in the different kernels have already been studied as part of $K_{112}^*$ or $K_{112}^o$.
So finally, combining the bound \eqref{bounddomega3} (corresponding in this case to the term $K_{112}^*$) and \eqref{boundmeanzero}, from \eqref{nabla2u3} we get the estimate
\begin{equation*}
\|\grad^2v_3\|_{L^1_T(L^\infty)}\leq c_1(\delta,\gamma,U)T+c_2(U)T^{3/2}+c_3(\delta,U)T^2.
\end{equation*}
\qed

\begin{rem}\label{remarkcurvature}
From the bounds above, we notice that we get indeed $v_3\in L^\infty(0,T;W^{2,\infty})$. Using \eqref{datoinicial} to interpolate, we find $\grad^2 v_1 \in L^q(0,T;W^{2,\infty})$, $q\in(1,2/(2-s))\subset(1,4/3)$. Taking into account \eqref{nolineal}, we finally prove that $u\in L^q(0,T;W^{2,\infty})$ for $q\in(1,2/(2-s))$.
\end{rem}
\medskip
\begin{rem}\label{remstriated}
We note now that the ideas of striated regularity  also show that the boundedness of the operator $\partial_i^\perp\partial_j(-\Delta)^{-1}\partial_k\partial_1(\partial_t-\Delta)^{-1}_0$ reduces to estimating $\Delta(\partial_t-\Delta^{-1})\te$. Indeed, as
\begin{equation*}
\widehat{\grad^2 v_3}(\xi,t)=\int_0^t \frac{\xi_i^\perp \xi_j \xi_k \xi_1}{|\xi|^2}e^{-(t-\tau)|\xi|^2}\hat\te(\xi,\tau)d\xi d\tau+ \frac12\begin{bmatrix}
      0 & -1  \\
      1 &  0  \\
     \end{bmatrix}\int_0^t \xi_k\xi_1e^{-(t-\tau)|\xi|^2}\hat\te(\xi,\tau)d\tau,
\end{equation*}
one can use the decomposition (Lemma 7.41, \cite{PDEBOOK})
\begin{equation*}
\xi_i\xi_j=a_{ij}(x)|\xi|^2 +\sum_{l}b_{ij}^l \xi_l (W(x)\cdot \xi),
\end{equation*}
where the function $a_{ij}$ are bounded, $b_{ij}$ are H\"older continuous and $W$ is a vector field along which $\te$ has some extra regularity (for example, the tangent vector field), so that $W\cdot\grad\te$ is H\"older continuous (in the case of a patch and the tangent vector field, $W\cdot \te\equiv 0$ in the sense of distributions). Thus one can write
\begin{equation*}
\widehat{\grad^2 u}=\int_0^t a_{ij}(x)\xi_k\xi_1 e^{-(t-\tau)|\xi|^2}\hat\te(\xi)d\xi+\int_0^t \sum_l b_{ij}^l\frac{\xi_k}{|\xi|^2}(W(x)\cdot\xi \hat\te (\xi))d \xi+ \frac12\begin{bmatrix}
      0 & -1  \\
      1 &  0  \\
     \end{bmatrix}\xi_k\xi_1e^{-t|\xi|^2}\te(\xi),
\end{equation*}
hence applying the above decomposition again, the estimation of $\grad^2 u$ in $L^1(0,T;L^\infty)$ reduces to estimate $\|\Delta(\partial_t-\Delta^{-1})\te\|_{L^1(0,T;L^\infty)}$ and to control the striated regularity of the second term, which can be done by choosing $W$ as the vector field tangent to the temperature patch (see Theorem 7.40 in \cite{PDEBOOK} or section \ref{Case33}).

The above ideas rely strongly on paradifferential calculus techniques. In addition, the quadratic form of the double Riesz transform is essential in the decomposition used, so we would still have to deal with $\Delta\heatinv$ (i.e., $\partial_k\omega_3$).
\end{rem}

\section{Persistence of $C^{2+\gamma}$ regularity}\label{Case33}

 To get the propagation of this extra regularity for the patch we cannot proceed as before: to use the trajectories, one would need the velocity to have $\gamma$-H\"older continuous second derivatives. However, it is sufficient to control this regularity in the tangential direction.

\begin{thm}\label{Case3}
Assume $\gamma\in(0,1)$. Let $D_0\subset \R^2$ be a bounded simply connected domain with boundary $\partial D_0\in C^{2+\gamma}$, $u_0\in H^{1+\gamma+s}$, $s\in(0,1-\gamma)$ a divergence-free vector field and $\te_0=1_{D_0}$. Then,
$$\te(x,t)=1_{D(t)}(x) \textup{ with } \partial D \in L^\infty(0,T;C^{2+\gamma}),$$
where $D(t)=X(D_0,t)$, and there exists a unique global solution $(u,\te)$ of \eqref{temperature}-\eqref{Boussinesq} such that
$$
u\in L^\infty(0,T;H^{1+\gamma+s})\cap L^2(0,T;H^{2+\mu})\cap L^p(0,T; H^{2+\delta})\cap L^q(0,T;W^{2, \infty}),
$$
for any $T>0$, $\mu<\frac12$, $\mu\leq \gamma+s$, $\delta<\frac12$, $1\leq p< 2/(1-(\gamma+s-\delta))$, $1\leq q< 2/(2-(\gamma+s))$.
\end{thm}

\medskip

Proof: As $u_0\in H^{1+\gamma+s}$ and $\theta_0\in H^{\delta}$, for any $\delta\in(0,1/2)$, from \cite{KUKA} one gets that $u\in L^\infty(0,T;H^{1+\mu})\cap L^2(0,T;H^{2+\mu})$ with $\mu<1/2$, $\mu\leq \gamma+s$. Using the splitting \eqref{vorsplitdef} to bootstrap, for the initial data we find that 
\begin{equation*}
\omega_0\in H^{\gamma+s}\Longrightarrow \omega_1\in L^\infty(0,T;H^{\gamma+s})\cap L^1(0,T;H^{2+\gamma+\tilde{s}}),\hspace{0.3cm}\tilde{s}\in(0,s),
\end{equation*}
while for the convection term it holds that
\begin{equation*}
\left.\begin{aligned}
u&\in L^\infty(0,T;H^{1+\mu})\\
\omega&\in L^\infty(0,T;H^\mu)\\
\omega&\in L^2(0,T;H^{1+\mu})
\end{aligned}\right\}\Longrightarrow 
\left.\begin{aligned}u\omega&\in L^\infty(0,T;H^\mu)\\
u\omega&\in L^2(0,T;H^{1+\mu})\end{aligned}\right\}
\Rightarrow 
\omega_2\in L^\infty(0,T;H^{1+\mu})\cap L^2(0,T;H^{2+\mu}).
\end{equation*}
The third term remains as before  $\omega_3\in L^\infty(0,T;H^{1+\delta})$. Hence, $\omega\in L^\infty(0,T;H^{\gamma+s})$ and therefore $u\in L^\infty(0,T;H^{1+\gamma+s})$. By interpolation,
 $v_1\in L^p(0,T;H^{2+\delta})$, where if $\gamma+s\geq1/2$ then $p=\frac{2}{1-(\gamma+s-\delta)}$ and if $\gamma+s\leq 1/2$ then $p=2\frac{1-(s-\tilde{s})}{1-(\gamma+s-\delta)-(s-\tilde{s})}$, so $u\in L^p(0,T;H^{2+\delta})$, where $1\leq p< \frac{2}{1-(\gamma+s-\delta)}$.

Let's consider the level-set characterization of the patch: 
\begin{equation*}
D_0=\{x\in\R^2: \varphi_0(x)>0\},
\end{equation*}
which for time $t$, $D(t)=X(t, D_0)$, is given by the function $\varphi(t,\cdot)$:
\begin{equation*}
\partial_t\varphi+u\cdot\grad\varphi=0,\hspace{0.7cm}
\varphi(x,0)=\varphi_0(x).
\end{equation*}
The vector field $W(t)=\grad^\perp \varphi(t)$ is then tangent to $\partial D(t)$ and its evolution is given by
\begin{equation}\label{tangent}
\begin{aligned}
\partial_t W +u\cdot \grad W=W\cdot\grad u,\hspace{0.5cm}
W(0)=\grad^\perp \varphi_0.
\end{aligned}
\end{equation}

In order to control the $C^{2+\gamma}$ regularity of $\partial D(t)$ one just need to show that $\grad W$ remains in $C^{\gamma}$.  By differentiating \eqref{tangent} one obtains

\begin{equation}\label{tangentder}
\partial_t \grad W + u\cdot\grad(\grad W) =W\cdot\grad(\grad u) + \grad W\cdot \grad u +\grad u \cdot \grad W.
\end{equation}
It is clear that we can choose $\varphi_0$ such that
\begin{equation}\label{campotangenteini}
W_0\in L^2\cap L^\infty,\hspace{0.4cm}\grad W_0 \in L^2\cap L^\infty.
\end{equation}
Then, from \eqref{tangent} and \eqref{tangentder} we deduce that  
\begin{equation}\label{campotangengel2}
W\in L^\infty(0,T;L^2\cap L^\infty), \hspace{0.3cm}\grad W\in L^\infty(0,T;L^2\cap L^\infty),
\end{equation}
since we know that $\grad u\in L^1(0,T;L^\infty)$, $\grad^2u\in L^1(0,T;L^\infty)$.

As $u$ is Lipschitz, the following estimate holds for all $t\in[0,T]$:
\begin{equation*}\label{striatedbound}
\begin{aligned}
\|\grad W\|_{C^\gamma}(t)&\!\leq\!\! \|\grad W_0\|_{C^\gamma}e^{\gamma\int_0^t\|\grad u\|_{L^\infty}d\tau}\!+\!e^{\gamma\int_0^t\|\grad u\|_{L^\infty}d\tau}\!\!\int_0^t\! \left( \|W\cdot\grad^2u\|_{C^\gamma}\!+\!2\|\grad W\|_{C^\gamma} \|\grad u\|_{C^\gamma}\!\right)d\tau.
\end{aligned}
\end{equation*}
From this and previous estimates we get that
\begin{equation}\label{holdergradW}
\|\grad W\|_{C^\gamma}(t)\leq c_1(T)+c_2(T)\int_0^t \|W\cdot\grad^2u\|_{C^\gamma}(\tau)d\tau+c_3(T)\int_0^t\|\grad u\|_{C^\gamma}(\tau)\|\grad W\|_{C^\gamma}(\tau)d\tau.
\end{equation}

We need to decompose the term $W\cdot\grad^2 u$ above to benefit from the extra cancellation in the tangential direction. We will use the following lemma: 
\medskip
\begin{lemma}\label{lemmafin}
For any $\gamma\in (0,1)$, there exists a constant $C$ such that the following estimate holds true:
\begin{equation}\label{striateddecomposition}
\begin{aligned}
\|W\cdot\grad^2 u\|_{L_t^1(C^\gamma)}&\leq C(\gamma)\left(\!\int_0^t \|W\|_{C^\gamma}\left( \|\grad\omega\|_{L^\infty}\!+\!\| \grad^2u\|_{L^\infty} \right) dt +\!\!\int_0^t\|\nabla\cdot (W\omega)\|_{C^\gamma}dt\right).
\end{aligned}
\end{equation}
\end{lemma}

Proof of Lemma \ref{lemmafin}: We decompose as follows
\begin{equation*}
\begin{aligned}
W(x)\cdot \grad(\partial_j u_i)(x)&=F(x)+\frac{1}{2\pi}\pv\int \frac{\sigma_{ij}(x-y)}{|x-y|^2}W(y)\cdot\grad  \omega(y) dy +\frac{W(x)\cdot\grad\omega(x)}{2}\begin{bmatrix}
      0 & -1  \\
      1 &  0  \\
     \end{bmatrix},
\end{aligned}
\end{equation*}
where
\begin{equation*}
F(x)=\frac{1}{2\pi}\pv\int \frac{\sigma_{ij}(x-y)}{|x-y|^2}(W(x)-W(y))\cdot \grad \omega(y) dy,\hspace{0.4cm} \sigma_{ij}(x)=\frac{1}{|x|^2} \begin{bmatrix}
2x_1x_2 & x_2^2-x_1^2  \\
x_2^2-x_1^2 &  -2x_1x_2  \\
\end{bmatrix}.
\end{equation*}
As the vector field $W$ is divergence-free, it yields that
\begin{equation*}
\|W\cdot\grad^2 u\|_{C^\gamma}\leq  \left|\left|F\right|\right|_{C^\gamma} + c \|\grad\cdot ( W\omega)\|_{C^\gamma}.
\end{equation*}
By the Lemma in Appendix of \cite{CONST}, 
\begin{equation*}
\|F\|_{C^\gamma} \leq C(\gamma)\|W\|_{C^\gamma}( \|\grad\omega\|_{L^\infty}+\| \grad^2u\|_{L^\infty}),
\end{equation*}
and the result follows.
\qed

\medskip
Hence, as $u\in L^1(0,T;W^{2,\infty})$ from Theorem \ref{Case2}, \eqref{striateddecomposition} yields
\begin{equation*}
\int_0^t\|W\cdot\grad^2u\|_{C^\gamma}(\tau)d\tau\leq c(T)+c(\gamma)\int_0^t\|\grad\cdot(W\omega)\|_{C^\gamma}d\tau.
\end{equation*}
Going back to \eqref{holdergradW},
\begin{equation}\label{secondbound}
\|\grad W\|_{C^\gamma}(t)\leq c_1(T)+c_2(T)\int_0^t\|\grad\cdot(W\omega)\|_{C^\gamma}(\tau)d\tau+c_3(T)\int_0^t\|\grad u\|_{C^\gamma}\|\grad W\|_{C^\gamma}(\tau)d\tau.
\end{equation}

\medskip
So the problem is reduced to control the tangential derivatives of the vorticity in $L^1(0,t;C^\gamma)$. We will use the properties of the heat kernel  in H\"older spaces applied to the equation satisfied by $\grad\cdot(W \omega)$:
\begin{equation*}
\begin{aligned}
\partial_t \grad\cdot(W\omega)+u\cdot\grad (\grad\cdot(W\omega))-\Delta \grad\cdot(W\omega) &= \grad\cdot (W\partial_1\theta)+\grad\cdot(W\Delta \omega-\Delta (W\omega)),\\
\grad\cdot(W\omega)|_{t=0}&=\grad\cdot(W_0\omega_0).
\end{aligned}
\end{equation*}
Notice that since $W$ is a divergence-free vector field tangent to the patch $\te$, it is possible to find
$$
\grad\cdot(W\partial_1\te)=\grad\cdot(\partial_1(W\te)-\te\partial_1W)=-\grad\cdot(\te\partial_1W).
$$
The solution can be decomposed as follows
\begin{equation}\label{decompositiondiv}
\begin{aligned}
\grad\cdot (W\omega) (t)&=e^{t\Delta}(\grad\cdot (W_0\omega_0)) -\heatinv\grad\cdot(u\grad\cdot (W\omega))\\
&\quad+\heatinv\grad\cdot (W\Delta \omega-\Delta (W\omega)) -\heatinv\grad\cdot(\te\partial_1W)\\
&=G_1+G_2+G_3+G_4.
\end{aligned}
\end{equation}
We will use now some classic estimates for solutions of the heat equation in H\"older spaces to achieve the gain of two derivatives integrable in time. We adapt these estimates to negative H\"older spaces (see Appendix). 
 We can apply \eqref{heatinitialneg} to $G_1$ to obtain
\begin{equation*}
\begin{aligned}
\|G_1\|_{L_t^1(C^\gamma)}&\leq c\left(\|\grad\cdot (W_0\omega_0)\|_{C^{-2+\gamma+s}}+\|W_0\omega_0\|_{L^1}\right)\leq c\left(\|W_0\omega_0\|_{C^{-1+\gamma+s}}+\|W_0\|_{L^2}\|\omega_0\|_{L^2}\right)\\
&=c\left(\|\grad^\perp \cdot(W_0\otimes u_0)-u_0\cdot \grad^\perp W_0\|_{C^{-1+\gamma+s}}+\|W_0\|_{L^2}\|\omega_0\|_{L^2}\right)\\
&\leq c\left( \|W_0\|_{C^{\gamma+s}}\|u_0\|_{C^{\gamma+s}}+\|u_0\|_{L^\infty}\|\grad^\perp W_0\|_{L^\infty}+\|W_0\|_{L^2}\|\omega_0\|_{L^2}\right),
\end{aligned}
\end{equation*}
thus 
\begin{equation}\label{boundG1}
\|G_1\|_{L_t^1(C^\gamma)}\leq c(t).
\end{equation}
Estimate \eqref{heatforcerneg} yields
\begin{equation*}
\begin{aligned}
\|G_2\|_{L_t^1(C^\gamma)}&\leq c\left(  \|\grad\cdot(u\grad\cdot(W\omega))\|_{L^r_t(C^{-2+\gamma})}+\|u\grad\cdot(W\omega)\|_{L^r_t(L^1)}\right)\\
&\leq c\left( \|\grad\cdot(u\otimes W\omega)-\omega W\cdot\grad u\|_{L_t^r(C^{-1+\gamma)}}+\| \|u\|_{L^2}\|W\|_{L^\infty}\|\grad\omega\|_{L^2}\|_{L^r_t}\right)\\
&\leq c\left( \|W\|_{L^\infty_t(C^\gamma)}\|u\|_{L^\infty_t(C^\gamma)}\|\omega\|_{L^r_t(C^\gamma)}+\|W\|_{L^\infty_t(L^\infty)}\|\omega\|_{L^{2r}_t(L^\infty)}\|\grad u\|_{L^{2r}_t(L^\infty)}\right.\\
&\quad+\left.\|W\|_{L^\infty_t(L^\infty)}\|u\|_{L^\infty_t(L^2)}\|\grad \omega\|_{L^r_t(L^2)}\right).
\end{aligned}
\end{equation*}
Now, as in Remark \ref{remarkcurvature}, since $v_1 \in L^1(0,T; H^{3+\gamma+\tilde{s}})$ and $v_1\in L^2(0,T;H^{2+\gamma+s})$, by interpolation we obtain $$u\in L^q(0,t;W^{2,\infty}),\hspace{0.6cm}\textup{with}\hspace{0.2cm}q\in \left[1, \frac{2}{2-(\gamma+s)}\right),$$
and similarly for any $0<\eps<\mu$
\begin{equation*}
\omega\in L^\infty(0,T;H^{\gamma+s})\cap L^2(0,T;H^{1+\mu})\Longrightarrow \omega \in L^{2r}(0,T;H^{1+\eps}),\hspace{0.2cm}1<r\leq \frac{1-\gamma-s+\mu}{1-\gamma-s+\eps}<1+\frac{\mu}{1-(\gamma+s)}.
\end{equation*}
Therefore, by choosing $1<r<\min\{2/(2-(\gamma+s)),1+\mu/(1-(\gamma+s))\}$ we get that
\begin{equation}\label{boundG2}
\|G_2\|_{L_t^1(C^\gamma)}\leq c(t).
\end{equation}
In the term $G_3$ we need to use the following commutator: $$W\Delta \omega-\Delta(W\omega)=-\omega\Delta W-2\grad W\cdot\grad\omega=-\grad\cdot(\omega\grad W)-\grad W\cdot\grad\omega,$$
\begin{equation*}
G_3=-\int_0^t e^{(t-\tau)\Delta}\grad\cdot(\grad\cdot(\omega\grad W))d\tau-\int_0^t e^{(t-\tau)\Delta}\grad\cdot(\grad W\cdot\grad\omega)d\tau = G_{31}+G_{32}.
\end{equation*}
Notice that in $G_{31}$ we find two derivatives applied to $\omega\grad W$. Therefore, proceeding as in \eqref{fourierside}, we can apply a similar estimate to \eqref{heatforcerneg} in order to find that
\begin{equation*}
\begin{aligned}
\|G_{31}\|_{L^1_t(C^\gamma)}&\leq c\left( \|\grad\cdot(\omega\grad W)\|_{L^r_t(C^{-1+\gamma})}+\|\mathcal{F}^{-1}(\chi(\xi)\widehat{\omega\grad W})\|_{L^r_t(L^\infty)}\right)\\
&\leq \left( \| \|\omega\|_{C^\gamma}\|\grad W\|_{C^\gamma}\|_{L^r_t}+\|\grad W\|_{L^\infty_t(L^2)}\|\omega\|_{L^r_t(L^2)}\right),
\end{aligned}
\end{equation*}
so  taking into account \eqref{campotangengel2}, we get
\begin{equation}\label{boundG31}
\|G_{31}\|_{L^1_t(C^\gamma)}\leq c \| \|\omega\|_{C^\gamma}\|\grad W\|_{C^\gamma}\|_{L^r_t}+c(t).\medskip
\end{equation}
For $G_{32}$ applying \eqref{heatforcesneg} yields that
\begin{equation*}
\begin{aligned}
\|G_{32}\|_{L^1_t(C^\gamma)}&\leq c\left( \|\grad W\cdot\grad\omega\|_{L^1_t(C^{-1+\gamma+s})}+\|\grad W\cdot\grad \omega\|_{L^1_t(L^1)}\right)\\
&\leq c\left( \|\grad W\|_{L^\infty(L^\infty)}\|\grad\omega\|_{L^1_t(L^\infty)}+\|\grad W\|_{L^\infty_t(L^\infty)}\|\grad\omega\|_{L^1_t(L^2)}\right),
\end{aligned}
\end{equation*}
and therefore 
\begin{equation}\label{boundG32}
\|G_{32}\|_{L^1_t(C^\gamma)}\leq c(t).
\end{equation}
Joining \eqref{boundG31} and \eqref{boundG32} we obtain
\begin{equation}\label{boundG3}
\|G_{3}\|_{L^1_t(C^\gamma)}\leq c(t)+c \| \|\omega\|_{C^\gamma}\|\grad W\|_{C^\gamma}\|_{L^r_t}.
\end{equation}
Finally, for the temperature term we use that $\grad\cdot(W\te)=W\cdot\grad \te\equiv 0$ in the sense of distributions:
\begin{equation*}
\|G_4\|_{L_t^1(C^\gamma)}\!\leq\! c\left(\! \|\!\grad\!\cdot\! (\te\partial_1W)\|_{L_t^1(C^{-2+\gamma+s})}\!+\!\|\te \partial_1W\|_{L^1_t(L^1)}\right)\!\leq\! c\|\grad W\|_{L^\infty_t(L^\infty)}\left(\!\|\theta\|_{L^1_t(L^\infty)}\!+\!\|\theta\|_{L^1_t(L^1)}\!\right),
\end{equation*}
so similarly to the others terms it is easy to find that
\begin{equation}\label{boundG4}
\|G_4\|_{L_t^1(C^\gamma)}\leq c(t).
\end{equation}

\medskip
Combining the above bounds \eqref{boundG1}, \eqref{boundG2}, \eqref{boundG3} and \eqref{boundG4}, it follows that
\begin{equation*}
\int_0^t \|\grad\cdot(W\omega)\|_{C^\gamma}(\tau)d\tau=\|G\|_{L^1_t(C^\gamma)}\leq c(t)+ c\| \|\omega\|_{C^\gamma}\|\grad W\|_{C^\gamma}\|_{L^r_t},
\end{equation*}
and going back to \eqref{secondbound} we conclude that
\begin{equation*}
\begin{aligned}
\|\grad W\|_{C^\gamma}(t)&\leq c(T)\!+\!c(T)\!\!\int_0^t\!\! \|\grad u\|_{C^\gamma}(\tau)\|\grad W\|_{C^\gamma}(\tau)d\tau\!+\!c(T)\left(\int_0^t\! \|\omega\|_{C^\gamma}^r(\tau)\|\grad W\|_{C^\gamma}^r(\tau)d\tau\right)^{1/r}\\
&\leq c(T)+c(T)\left(\int_0^t \|\grad u\|_{C^\gamma}^r(\tau)\|\grad W\|_{C^\gamma}^r(\tau)d\tau\right)^{1/r}\\
&\leq c(T)\!+\!c(T)\left(\int_0^t\!\|\grad u\|_{C^\gamma}^{r_1}\right)^{1/r_1}\!\left(\int_0^t\|\grad W\|_{C^\gamma}^{r_2}\right)^{1/r_2}\!\!\leq c(T)\!+\!c(T)\left(\int_0^t\!\|\grad W\|_{C^\gamma}^{r_2}\!\right)^{1/r_2}\!\!,
\end{aligned}
\end{equation*}
where we can choose $1<r<r_1<\min\{2/(2-(\gamma+s)),1+\mu/(1-(\gamma+s))\}$ and $1/r_2=1/r-1/r_1$.
Therefore, raising to power $r_2$, 
\begin{equation*}
\|\grad W\|_{C^\gamma}^{r_2}(t)\leq c(T)+c(T)\int_0^t\|\grad W\|_{C^\gamma}^{r_2}(\tau)d\tau,
\end{equation*}
 and hence applying Gr\"onwall's inequality the result follows.
\qed

\begin{rem} The result also holds for $s=0$, i.e., $u_0\in H^{1+\gamma}$, using paradifferential calculus and the {\em{tilde spaces}} introduced in \cite{CHEMINTILDE}. In fact, Proposition \ref{propo3} in Appendix applied to \eqref{tangentder} yields that
\begin{equation}
\|\grad W\|_{L_T^\infty(C^\gamma)}\leq c(T)+c(T)\|W\cdot\grad^2 u\|_{\tilde{L}^1_t(C^\gamma)}+c(T)\int_0^t\|\grad W\|_{C^\gamma}\|\grad u\|_{C^\gamma}d\tau.
\end{equation}
Lemma \ref{lemmafin} can be written using $\tilde{L}^1_T$ instead of $L_T^1$, so applying Proposition \ref{propoheat} to the term $\grad\cdot(W\omega)$ yields the result for $s=0$.
\end{rem}

\section{Corollaries}\label{Corollaries}

The results above have been presented using only Sobolev spaces for the velocity, but as commented before they can be stated using also H\"older or Besov spaces. 
\begin{cor}\label{cor1}
Let $D_0\subset\R^2$ be a bounded simply connected domain with boundary $\partial D_0\in W^{2,\infty}$,  $u_0 \in C^{s_1}\cap L^p$ a divergence-free vector field, $s_1\in(0,1)$, $p\in[1,2)$ and $\te_0=1_{D_0}$ the characteristic function of $D_0$. Then, 
 $$\theta(x,t)=1_{D(t)}(x) \textup{ and }\partial D\in L^\infty(0,T;W^{2,\infty}).$$
Moreover, there exists a unique global solution $(u, \te)$ of \eqref{temperature}-\eqref{Boussinesq} such that $$u\in L^\infty(0,T;H^{s_2})\cap L^2(0,T;H^{1+s_2})\cap L^1(0,T;H^{2+s_2})\cap L^\infty(0,T;C^{s_1})\cap L^1(0,T;W^{2,\infty}),$$
for any $T>0$, $0\leq s_2<s_1/2<1/2$.
\end{cor}

Proof: First we notice that as $u_0\in C^{s_1}\cap L^p$, it also holds that $u_0\in L^{\tilde{p}}\cap H^{s_2}$ for $\tilde{p}\in[p, \infty]$, $s_2\in [0,s_1/2)$.
The {\em{a priori}} estimates are done exactly as in Theorem \ref{Case1}:
 $$u\in L^\infty(0,T;H^{s_2})\cap L^2(0,T;H^{1+s_2})\cap L^1(0,T; H^{2+s_2}).$$ 
In addition, as $u_0\in C^{s_1}$, $s_1>0$, one trivially gets  $\grad^2 v_1\in L^1(0,T;L^\infty)$.
 To show that $u\in L^1(0,T;W^{2,\infty})$, we follow the same structure as in Theorem \ref{Case2}: use the splitting \eqref{vorsplitdef} and treat each one of the terms $\grad^2 v_1$, $\grad^2v_2$ and $\grad^2v_3$, given by  \eqref{hesu} separately. 

   To deal with the temperature term, it suffices to prove that $u\in L^\infty(0,T;L^\infty)$ to apply the steps of Theorem \ref{Case2} to obtain $\grad^2 v_3\in L^\infty(0,T;L^\infty)$. Proposition 2.1 in \cite{DANCHINZHANG} would yield that $u\in L^\infty(0,T; C^{s_1})$. We give a different proof that does not require the use of paradifferential calculus. We apply the Leray projector to the velocity equation to obtain
  \begin{equation*}
  u_t-\Delta u=-\grad\cdot(u\otimes u)+\grad(-\Delta)^{-1}\left(\grad\cdot\grad\cdot(u\otimes u)\right)+(0,\theta)-\grad(-\Delta)^{-1}\partial_2\theta.
  \end{equation*}
  Hence, we can write the velocity as follows:
  \begin{equation}\label{udecomp}
  \begin{aligned}
    u=w_1+w_2+w_3+w_4+w_5,\hspace{5cm}\\
  w_1=e^{t\Delta}u_0,\hspace{0.4cm}w_2=-(\partial_t-\Delta)^{-1}_0\grad\cdot(u\otimes u),
  \hspace{0.4
  	cm} w_3=(\partial_t-\Delta)^{-1}_0\grad (-\Delta)^{-1}\grad\cdot(\grad\cdot(u\otimes u)), \\
  w_4=(\partial_t-\Delta)^{-1}_0(0,\theta),\hspace{0.5cm}w_5=-(\partial_t-\Delta)^{-1}_0\grad(-\Delta)^{-1}\partial_2\theta.\hspace{1.8cm}
  \end{aligned}
  \end{equation}
  From the properties of the heat equation we obtain
  \begin{equation}\label{boundw1w4}
  \|w_1\|_{L^\infty_T(C^{s_1})}\leq c\|u_0\|_{C^{s_1}},\hspace{0.4cm}\|w_4\|_{L^\infty_T(B_{\infty,\infty}^{2})}\leq c\|\theta\|_{L^\infty_T(L^\infty)},\hspace{0.4cm} \|w_5\|_{L^\infty_T(B_{\infty,\infty}^{2})}\leq c\|\theta\|_{L^\infty_T(L^\infty)}.
  \end{equation}
 Using the boundedness of singular integrals in H\"older spaces, $w_3$ can be treated as $w_2$. First, as $u\in L^\infty(0,T;H^{s_2})\cap L^2(0,T;H^{1+s_2})$, by interpolation $$u\in L^q(0,T;H^{1+\tilde{s}_2})\hookrightarrow L^q(0,T;L^\infty),\hspace{0.5cm}  q=\frac{2}{1-(s_2-\tilde{s}_2)}.$$   
 Now we proceed as follows
  \begin{equation}\label{truco}
  \|w_2\|_{C^{s_1}}(t)\leq c\int_0^t\frac{\|u\otimes u\|_{C^{s_1}}(\tau)}{(t-\tau)^{1/2}}d\tau\leq c\int_0^t\frac{\|u\|_{C^{s_1}}\|u\|_{L^\infty}}{(t-\tau)^{1/2}}d\tau.
  \end{equation}
Choose $l=2/(1+\eps)$ with $\eps\in(0,1)$ and $$\frac{1}{r}=1-\frac{1}{q}-\frac{1}{l}=\frac{s_2-\tilde{s}_2-\eps}{2},$$
so applying H\"older inequality yields that
\begin{equation}\label{boundw2}
\|w_2\|_{C^{s_1}}(t)\leq c(T) \|u\|_{L^q_T(L^\infty)}\left(\int_0^t \|u\|^r_{C^{s_1}}(\tau)d\tau \right)^{1/r}.
\end{equation}
From the decomposition \eqref{udecomp}, the bounds \eqref{boundw1w4}, \eqref{boundw2} and recalling that $\|w_3\|_{C^{s_1}}\leq c\|w_2\|_{C^{s_1}}$, it is easy to get
\begin{equation*}
\|u\|_{C^{s_1}}(t)\leq c(T)+c(T)\left(\int_0^t \|u\|^r_{C^{s_1}}(\tau)d\tau \right)^{1/r}.
\end{equation*}
Raising to the power $r$ and applying Gr\"onwall's inequality we conclude that
\begin{equation}\label{corfin}
u\in L^\infty(0,T;C^{s_1})
\end{equation}
and therefore $u\in L^\infty(0,T;L^\infty)$.

To conclude, we need to prove that $\grad^2 v_2 \in L^1(0,T;L^\infty)$. It suffices to show that $u\omega\in L^r(0,T;C^{\delta})$ for some $r>1,\delta>0$. By interpolation,
\begin{equation}
\omega\in L^\sigma(0,T; H^{1+s_4}),\hspace{0.2cm}\sigma=\frac{2}{1+s_4}>\frac43, \hspace{0.2cm} s_4\in(0,s_2).
\end{equation}
Thus $\omega \in L^\sigma (0,T; C^{s_4})$. From this and \eqref{corfin}, $u\omega\in L^\sigma(0,T; C^\delta)$, with $\delta=\min\{s_3,s_4\}>0$.

\qed

\begin{cor}
Let $D_0\subset\R^2$ be a bounded simply connected domain with boundary $\partial D_0\in W^{2,\infty}$, $\te_0=1_{D_0}$ and $u_0$ a divergence-free vector field in the Besov space $B^1_{2,1}$. Then, 
$$\partial D\in L^\infty(0,T;W^{2,\infty})$$
and there exists a unique global solution $(u,\te)$ of \eqref{temperature}-\eqref{Boussinesq} such that
$$
u\in L^\infty(0,T;B^1_{2,1})\cap L^2(0,T;B^2_{2,1})\cap L^1(0,T; B^{2+\mu}_{2,1})\cap L^1(0,T;W^{2,\infty}),
$$
for any $T>0$, $\mu<1/2$.
\end{cor}

Proof: One can repeat the proof of Theorem \ref{Case2} to get this result. In this case, the {\em{a priori}} estimates are done using the result in \cite{DANCHIN}. Then, by splitting the vorticity equation, one can take advantage of the fact that $B^1_{2,1}\hookrightarrow L^\infty$, it is an algebra and for the heat equation with initial data in this space there is a gain of two derivatives integrable in time \cite{FERO3}. \qed

\begin{cor}\label{cor3}
Let $D_0\subset\R^2$ be a bounded simply connected domain with boundary $\partial D_0\in C^{1+\gamma}$, $u_0\in B^{-1+\gamma+s_1}_{\infty,\infty}\cap H^{s_2}$ a divergence-free vector field, $s_1\in(0,1-\gamma),s_2\in(0,1)$ and $\theta_0=1_{D_0}$ the characteristic function of $D_0$. Then, there exists a unique global solution $(u, \te)$ of \eqref{temperature}-\eqref{Boussinesq} such that	
 $$\theta(x,t)=1_{D(t)}(x) \textup{ and }\partial D\in L^\infty(0,T;C^{1+\gamma} ).$$
 Moreover,
 $$u\in L^\infty(0,T;H^{s_2})\cap L^2(0,T;H^{1+s_2})\cap L^1(0,T;H^{2+\mu})\cap L^\infty(0,T;B_{\infty,\infty}^{-1+\gamma+s_1})\cap L^1(0,T;C^{1+\gamma+\tilde{s}_1}),$$
 for any $T>0$, $\mu<\min\{\frac12,s_2\}$, $0<\tilde{s}_1<s_1$.
\end{cor}

Proof: Since $u_0\in H^{s_2}$, we get the {\em{a priori}} estimates $u\in L^\infty(0,T;H^{s_2})\cap L^2(0,T;H^{1+s_2})\cap L^1(0,T;H^{2+\mu})$. Now, we use the decomposition \eqref{udecomp}. As $u_0\in B_{\infty,\infty}^{-1+\gamma+s_1}$,  it holds that $w_1 \in L^\infty(0,T;C^{-1+\gamma+s_1})\cap L^2(0,T;C^{\gamma+\sigma})\cap L^1(0,T;C^{1+\gamma+\tilde{s}_1})$, $\tilde{s}_1<\sigma<s_1$. Since $w_4, w_5\in L^\infty(B_{\infty,\infty}^2)$, we only need to deal with the nonlinear terms. 

As $w_2=\heatinv \grad\cdot (u\otimes u)$, it suffices to show that $u\in L^2(0,T;C^{\gamma+\sigma})$ to conclude that $w_2\in L^1(0,T;C^{1+\gamma+\tilde{s}_1})$, $\tilde{s}_1\in(0,s_1)$. By the estimates above, we only need to show that $w_2\in L^2(0,T;C^{\gamma+\sigma})$. This means in turn that it suffices to show $u\otimes u\in L^2(0,T;B_{\infty,\infty}^{-1+\gamma+\sigma})$. We will show now that $u\in L^\infty(0,T; B^{-1+\gamma+s_1})$ indeed. It is clear for $w_1, w_4,w_5$. 

Applying basic paradifferential calculus estimates (see Chapter $2$ in \cite{PDEBOOK}), we obtain
\begin{equation}\label{paracal}
\|u\otimes u\|_{B_{\infty,\infty}^{-1+\gamma+s_1}}\leq c \|u\|_{B_{\infty,\infty}^{-1+\gamma+s_1}}(\tau)(\|u\|_{L^\infty}(\tau)+\|u\|_{H^1}).
\end{equation}
Proceeding now as in \eqref{truco},
\begin{equation*}
\|w_2\|_{B_{\infty,\infty}^{-1+\gamma+s_1}}(t)\leq c\int_0^t\frac{\|u\otimes u\|_{B_{\infty,\infty}^{-1+\gamma+s_1}}(\tau)}{(t-\tau)^{1/2}}d\tau\leq c\int_0^t \frac{\|u\|_{B_{\infty,\infty}^{-1+\gamma+s_1}}(\tau)(\|u\|_{L^\infty}(\tau)+\|u\|_{H^1}(\tau))}{(t-\tau)^{1/2}}d\tau,
\end{equation*}
we find that (for the same $r$ as in \eqref{boundw2})
\begin{equation*}
\|u\|_{B_{\infty,\infty}^{-1+\gamma+s_1}}(t)\leq c(T)+c(T)\left(\int_0^t \|u\|^r_{B_{\infty,\infty}^{-1+\gamma+s_1}}(\tau)d\tau\right)^{1/r},
\end{equation*}
so we conclude that $u\in L^\infty(0,T;B_{\infty,\infty}^{-1+\gamma+s_1})$. From this and the inequality \eqref{paracal} we conclude that $u\otimes u\in L^2(0,T;B^{-1+\gamma+s_1}_{\infty,\infty})$.

\qed

\vspace{0.2cm}

 \section*{\hfil Appendix\hfil}

We include here some results that have been used related to the linear transport and heat equations in Sobolev and H\"older spaces. At the end we include the kernels introduced in \eqref{kernelsref} to define the operators $R_iR_j\partial_k\partial_1\heatinv$. 

First, we recall the definition of the Besov spaces (see \cite{PDEBOOK} for details). Consider the nonhomogeneous Littlewood-Paley decomposition:
\newline
\newline
 Let $B=\{|\xi|\in\R^2: |\xi|\leq 4/3\}$ and $C=\{|\xi|\in\R^2:3/4\leq|\xi|\leq8/3\}$, and fix two smooth radial functions $\chi$ and $\varphi$ supported in $B$, $C$, respectively, and satisfying 
\begin{equation}\label{littlewoodp}
\chi(\xi)+\sum_{j\geq0}\varphi(2^{-j}\xi)=1,\hspace{0.3cm}\forall \xi\in \R^2.
\end{equation}
The nonhomogeneous dyadic blocks are defined as $\Delta_jf=\mathcal{F}^{-1}(\varphi(2^{-j}\xi)\hat{f}(\xi))$ for $j\geq0$ and $\Delta_{-1}f=\mathcal{F}^{-1}(\chi(\xi)\hat{f}(\xi))$.
Then, the Besov space $B_{p,q}^\gamma(\R^2)$, $\gamma\in\R$, $p,q\in [1,\infty]$ is defined by
$$B_{p,q}^\gamma(\R^2)=\{u\in S'(\R^2): \|u\|_{B_{p,q}^\gamma}=\|2^{j\gamma}\|\Delta_j u\|_{L^p}\|_{l^q(j\geq-1)}<\infty\},$$
where $S'(\R^2)$ denotes the space of tempered distributions over $\R^2$. We recall that $H^s=B_{2,2}^s$ and $C^\gamma=B_{\infty,\infty}^\gamma$ for $s\in \R$, $\gamma\in(0,1)$.

\medskip
\begin{prop}\label{propoinitial} Let $s>0$, $r\in(1,\infty]$. Then, the following estimates hold
	\begin{equation}\label{heatforcer}
	\|\partial_i\partial_k(\partial_t-\Delta)^{-1}_0f\|_{L^r_T(C^\gamma)}\leq c\|f\|_{L^r_T(C^{\gamma})},
	\end{equation}
	\begin{equation}\label{heatforces}	
	\|\partial_i\partial_k(\partial_t-\Delta)^{-1}_0f\|_{L^1_T(C^\gamma)}\leq c\|f\|_{L^1_T(C^{\gamma+s})},
	\end{equation} 
	\begin{equation}\label{heatinitial}
	\|\partial_i\partial_ke^{t\Delta}u_0\|_{L^1_T(C^\gamma)}\leq c \|u_0\|_{C^{\gamma+s}}.
	\end{equation}
	Furthermore, there exists $u_0\in C^\gamma$ for which $\partial_i\partial_ke^{t\Delta}u_0\notin L^1_T(C^\gamma).$
	\end{prop}

Proof: The proof of \eqref{heatforcer} can be found in \cite{KRYLOV}.  The proof \eqref{heatforces}  follows from Bernstein inequalities and the decay of the heat kernel:
\begin{equation*}
\begin{aligned}
\|\partial_i\partial_k\heatinv f\|_{L^1_T(C^{\gamma})}&\leq c\int_0^T \sup_{j\geq -1}2^{j(\gamma+s)}2^{j(2-s)}\int_0^t e^{-c(t-\tau)2^{2j}}\|\Delta_j f\|_{L^\infty}d\tau dt\\
&\leq c\int_0^T \int_0^t \frac{c}{(t-\tau)^{1-s/2}}\sup_{j\geq-1}2^{j(\gamma+s)}\|\Delta_j f\|_{L^\infty}d\tau dt\leq c(T)\|f\|_{L^1_T(B_{\infty,\infty}^{\gamma+s})}.
\end{aligned}
\end{equation*}
We get \eqref{heatinitial} as before
\begin{equation*}
\begin{aligned}
\|\partial_i\partial_ke^{t\Delta}u_0\|_{L^1_T(C^\gamma)}&\leq \!c\!\! \int_0^T\!\!\!\!\sup_{j\geq -1}2^{j(\gamma+s)}2^{j(2-s)}e^{-ct2^{2j}}\!\|\Delta_j u_0\|_{L^\infty}dt\leq  \!c\!\int_0^T \!\frac{\|u_0\|_{C^{\gamma+s}}}{t^{1-s/2}}dt\leq\! c(T)\|u_0\|_{C^{\gamma+s}}.
\end{aligned}
\end{equation*}
To prove the last statement we proceed as follows
\begin{equation*}
\begin{aligned}
\|\partial_i\partial_ke^{t\Delta}u_0\|_{C^\gamma}\geq c \sup_{j\geq 0}2^{j\gamma}2^{2j}\|\Delta_j (e^{t\Delta}u_0)\|_{L^\infty}\geq c\sup_{j\geq 0}2^{j\gamma}2^{2j}e^{-ct2^{2j}}\|\Delta_ju_0\|_{L^\infty},
\end{aligned}
\end{equation*}
therefore
\begin{equation*}
\begin{aligned}
\|\partial_i\partial_ke^{t\Delta}\|_{L^1_T(C^\gamma)}\geq c\left|\left|\sup_{j\geq 0}\left(2^{2j}e^{-ct2^{2j}}\right)\inf_{j\geq0}\left(2^{j\gamma}\|\Delta_ju_0\|_{L^\infty}\right)\right|\right|_{L^1_T}=c\left|\left|\frac{1}{t}\inf_{j\geq0}\left(2^{j\gamma}\|\Delta_ju_0\|_{L^\infty}\right)\right|\right|_{L^1_T}.
\end{aligned}
\end{equation*}
Thus one only needs to find $u_0\in C^\gamma$ such that $\inf_{j\geq0}\left(2^{j\gamma}\|\Delta_ju_0\|_{L^\infty}\right)>0$. It is not difficult to check that the function define by $\displaystyle\hat{u}_0(\xi)=\frac{1}{|\xi|^{2+\gamma}}\left(1-\chi(\xi)\right)$ satisfies the condition.
\qed

\vspace{0.3cm}

Now we adapt these estimates to negative H\"older spaces.

\begin{prop} 
Let $g(x,t)=\grad\cdot f(x,t)$, then we have the following estimates
	\begin{equation}\label{heatforcerneg}
	\left|\left| (\partial_t-\Delta)^{-1}_0 g \right|\right|_{L^r_T(C^\gamma)} \leq c \left( \|g\|_{L^r_T(B^{-2+\gamma}_{\infty,\infty})} +\|f\|_{L^r_T(L^1)}\right),
	\end{equation}
	\begin{equation}\label{heatforcesneg}
	\left|\left| (\partial_t-\Delta)^{-1}_0 g \right|\right|_{L^1_T(C^\gamma)}\leq c \left( \|g\|_{L^1_T(B^{-2+\gamma+s}_{\infty,\infty})} +\|f\|_{L^1_T(L^1)}\right),
	\end{equation}
	\begin{equation}\label{heatinitialneg}
	\left|\left|e^{t\Delta} g \right|\right|_{L^1_T(C^\gamma)}\leq c \left( \|g\|_{B^{-2+\gamma+s}_{\infty,\infty}} +\|f\|_{L^1}\right).
	\end{equation}
\end{prop}

Proof: First, if we denote by $\chi$ the first dyadic block in the Littlewood-Paley decomposition \eqref{littlewoodp}, we notice that 
$$g\in B^{-2+\gamma}_{\infty,\infty}\Leftrightarrow \exists h\in C^\gamma\hspace{0.2cm} {\rm{such \hspace{0.1cm}that}}\hspace{0.2cm} g=\Delta h.$$
Indeed,
\begin{equation}\label{fourierside}
\begin{aligned}
\|h\|_{C^\gamma}&=2^{-\gamma}\left|\left|\mathcal{F}^{-1}\left(\chi(\xi)\frac{\xi}{|\xi|^2}\cdot \hat{f}(\xi)\right)\right|\right|_{L^\infty}+\sup_{j\geq 0}2^{j\gamma}\|\Delta_j (\Delta^{-1}g)\|_{L^\infty} \\
&\leq2^{-\gamma}\left|\left|\mathcal{F}^{-1}\left(\chi(\xi)\frac{\xi}{|\xi|^2}\cdot \hat{f}(\xi)\right)\right|\right|_{L^\infty}+ \|g\|_{B^{-2+\gamma}_{\infty,\infty}}.
\end{aligned}
\end{equation}
Note that if $f\in L^1$,
\begin{equation*}
\begin{aligned}
\left|\left|\mathcal{F}^{-1}\left(\chi(\xi)\frac{\xi}{|\xi|^2}\cdot \hat{f}(\xi)\right)\right|\right|_{L^\infty} \leq c \|f\|_{L^1},
\end{aligned}
\end{equation*}
so it holds that
\begin{equation*}
\|h\|_{C^\gamma} \leq  c\|f\|_{L^1}+\|g\|_{B^{-2+\gamma}_{\infty,\infty}}.
\end{equation*}

\medskip
Applying the classic estimates  \eqref{heatforcer}, \eqref{heatforces} and \eqref{heatinitial} to $\heatinv g=\Delta \heatinv h$ and $e^{t\Delta}g=\Delta e^{t\Delta}h$, \eqref{heatforcerneg}, \eqref{heatforcesneg} and \eqref{heatinitialneg} follow.
\qed

\medskip
Now we show two results that allow us to get the $C^{2+\gamma}$ persistence of regularity for $u_0\in H^{1+\gamma}$. We need to recall the definition of the {\em{tilde spaces}}:
$$\tilde{L}^\rho_t(B_{p,q}^\gamma(\R^2))=\{u \in S'(0,t;\R^2): \|u\|_{\tilde{L}^\rho_t(B_{p,q}^\gamma)}=\|2^{j\gamma}\|\Delta_j u\|_{L^\rho_t(L^p)}\|_{l^q(j\geq-1)}<\infty\}.$$

\begin{prop}\label{propoheat} Let $f\in \tilde{L}^1_T(C^\gamma)$, $f_0\in C^\gamma$. Then,
\begin{equation*}
\begin{aligned}
\|\partial_i\partial_k(\partial_t-\Delta)_0^{-1}f\|_{\tilde{L}_T^1(C^\gamma)}&\leq c\|f\|_{\tilde{L}_T^1(C^\gamma)}\leq c\|f\|_{L_T^1(C^\gamma)},\\
\|\partial_i\partial_k e^{t\Delta}f_0\|_{\tilde{L}^1_T(C^\gamma)}&\leq c(T)\|f_0\|_{C^\gamma}.
\end{aligned}
\end{equation*}
\end{prop}

Proof: By the definition of {\em{tilde}} spaces we can now integrate first in time, 
\begin{equation*}
\begin{aligned}
\|\partial_i\partial_k(\partial_t&-\Delta)_0^{-1}f\|_{\tilde{L}_T^1(C^\gamma)}=\sup_{j\geq-1}2^{j\gamma}\|\Delta_j(\partial_i\partial_k(\partial_t-\Delta)_0^{-1}f)\|_{L_T^1(L^\infty)}\\
&\leq\! c\sup_{j\geq-1} \left|\left|\!\int_0^t\! 2^{2j}e^{-c(t-\tau)2^{2j}}2^{j\gamma}\|\Delta_j f\|_{L^\infty}(\tau) d\tau\right|\right|_{L^1_T}\!\!\leq\! c\sup_{j\geq -1}2^{j\gamma}\|\Delta_j f\|_{L^1_T(L^\infty)}=c\|f\|_{\tilde{L}^1_T(C^\gamma)}.
\end{aligned}
\end{equation*}
The initial condition estimate follows from a similar procedure.
\qed
\begin{prop}\label{propo3} Let $f$ be the solution to the transport equation
$$f_t+u\cdot\grad f=g, \hspace{0.3cm}f(0)=f_0,$$
where $\grad u\in L^1_T(L^\infty)$, $g\in \tilde{L}^1_T(C^\gamma)$ and $f_0\in C^\gamma$. Then,
$$\|f\|_{L^\infty_T(C^\gamma)}\leq\left(\|f_0\|_{C^\gamma}+\|g\|_{\tilde{L}_T^1(C^\gamma)}\right)e^{c\int_0^t}\|\grad u\|_{L^\infty}d\tau.$$
\end{prop}
The proof of Proposition \ref{propo3} can be found in Theorem $3.14$, Chapter 3 \cite{PDEBOOK}. 
\qed

\vspace{0.5cm}

We include next the explicit expression of the kernels $K_{ijk}$ \eqref{kernelsref}:

\begin{equation}\label{allthekernels}
\begin{aligned}
K_{111}&=\left(\frac{24x_1^3x_2}{\pi |x|^5}-\frac{12x_1x_2}{\pi|x|^3}\right)G(|x|,t)-\frac{e^{-|x|^2/4t}}{\pi(4t)^2}\left(\frac{12x_1^3x_2}{|x|^4}-\frac{6x_1x_2}{|x|^2}+\frac{4x_1^3x_2}{|x|^2}\frac{1}{4t}\right),\\
K_{112}&=\left(\frac{24x_1^2x_2^2}{\pi |x|^5}-\frac{3}{\pi|x|}\right)G(|x|,t)-\frac{e^{-|x|^2/4t}}{\pi(4t)^2}\left(-2+\frac{4}{4t}\frac{x_1^2x_2^2}{|x|^2}+12\frac{x_1^2x_2^2}{|x|^4}\right),\\
K_{122}&=\left(-\frac{12x_1x_2}{\pi|x|^3}+\frac{24x_1x_2^3}{\pi|x|^5}\right)G(|x|,t)-\frac{e^{-|x|^2/4t}}{\pi(4t)^2}\left(-\frac{12x_1x_2}{|x|^2}+\frac{12x_1x_2^3}{|x|^4}+\frac{4x_1x_2^3}{|x|^2}\frac{1}{4t}\right),\\
K_{211}&=\left(-\frac{24x_1^4}{\pi|x|^5}+\frac{24x_1^2}{\pi|x|^3}-\frac{3}{\pi|x|}\right)G(|x|,t)-\frac{e^{-|x|^2/4t}}{\pi(4t)^2}\left(\frac{12x_1^2}{|x|^2}-\frac{12x_1^4}{|x|^4}-\frac{4x_1^4}{|x|^2}\frac{1}{4t}\right),\\
K_{121}&=K_{112}, \hspace{0.4cm}K_{212}=-K_{111},\hspace{0.4cm} K_{221}=-K_{111},\hspace{0.4cm} K_{222}=-K_{112},
\end{aligned}
\end{equation}
\medskip
where we recall that
$$G(|x|,t)=\frac{1}{|x|^3}(1-e^{-|x|^2/4t})-\frac{1}{4t |x|}e^{-|x|^2/4t}\geq 0\hspace{0.5cm}\forall |x|,t\geq 0.$$

\subsection*{{\bf Acknowledgments}}
This research was partially supported by the project P12-FQM-2466 of Junta de Andaluc\'ia, Spain, 
 the grant MTM2014-59488-P (Spain) and by the ERC through the Starting Grant project H2020-EU.1.1.-639227. 
EGJ was supported by MECD FPU grant from the Spanish Government. FG acknowledges support from the Ram\'on y Cajal program RyC-2010-07094.

\vspace{1cm}

\begin{quote}
\begin{tabular}{ll}
\textbf{Francisco Gancedo}\\
{\small Departamento de An\'{a}lisis Matem\'{a}tico $\&$ IMUS}\\
{\small Universidad de Sevilla}\\
{\small C/ Tarfia s/n, Campus Reina Mercedes, 41012 Sevilla, Spain}\\
{\small Email: fgancedo@us.es}
\end{tabular}
\end{quote}

\begin{quote}
\begin{tabular}{ll}
\textbf{Eduardo Garc\'ia-Ju\'arez}\\
{\small Departamento de An\'{a}lisis Matem\'{a}tico $\&$ IMUS}\\
{\small Universidad de Sevilla}\\
{\small C/ Tarfia s/n, Campus Reina Mercedes, 41012 Sevilla, Spain}\\
{\small Email: eduardogarcia@us.es}
\end{tabular}
\end{quote}

\end{document}